\makeatletter
\def\cl@chapter{}
\makeatother

\documentclass[epjST]{svjour}

\newcommand\figurespath{figures/}
\usepackage{graphicx}
\usepackage{svg}
\graphicspath{{\figurespath}}

\usepackage[english]{babel}
\usepackage[utf8]{inputenc}

\usepackage{calc}
\usepackage{amsmath}
\usepackage{amssymb}
\usepackage{float}
\usepackage{subfig}
\usepackage[export]{adjustbox}
\usepackage{hyperref}
\usepackage{csquotes}
\floatstyle{ruled}
\newfloat{algorithm}{tbp}{loa}
\providecommand{\algorithmname}{Algorithm}
\floatname{algorithm}{\protect\algorithmname}

\usepackage[backend=biber,style=phys,articletitle=false,biblabel=brackets,chaptertitle=false,pageranges=false, natbib=true, doi=true, bibencoding=utf8]{biblatex}
\usepackage{amsmath,amssymb}
\usepackage[bb=boondox]{mathalfa} 
\DeclareFontFamily{U}{BOONDOX-calo}{\skewchar\font=45 }
\DeclareFontShape{U}{BOONDOX-calo}{m}{n}{
  <-> s*[1.05] BOONDOX-r-calo}{}
\DeclareFontShape{U}{BOONDOX-calo}{b}{n}{
  <-> s*[1.05] BOONDOX-b-calo}{}
\DeclareMathAlphabet{\mathcalboondox}{U}{BOONDOX-calo}{m}{n}
\SetMathAlphabet{\mathcalboondox}{bold}{U}{BOONDOX-calo}{b}{n}
\DeclareMathAlphabet{\mathbcalboondox}{U}{BOONDOX-calo}{b}{n}

\usepackage{cleveref}
\Crefname{equation}{Eq.}{Eqs.}
\crefname{equation}{Eq.}{Eqs.}
\Crefname{figure}{Fig.}{Figs.}
\crefname{figure}{Fig.}{Figs.}
\Crefname{tabular}{Tab.}{Tabs.}
\crefname{tabular}{Tab.}{Tabs.}
\Crefname{appendix}{Appendix}{Appendices}
\crefname{appendix}{Appendix}{Appendices}

\newcommand\tsm{\textsc{tsm}}

\renewcommand\lg{\textsc{lg}}
\newcommand\et{\textsc{et}}

\newcommand{\A}{\ensuremath{\mathcalboondox{A}}}
\newcommand{\B}{\ensuremath{\mathcalboondox{B}}}
\newcommand{\x}{{\ensuremath{\mathcalboondox{x}}}}
\newcommand{\y}{{\ensuremath{\mathcalboondox{y}}}}
\newcommand{\X}{{\ensuremath{\mathcalboondox{X}}}}
\newcommand{\Xp}{{\ensuremath{\X^+}}}
\newcommand{\Xm}{{\ensuremath{\X^-}}}
\newcommand{\Xpm}{{\ensuremath{\X^{+/-}}}}
\newcommand{\Y}{{\ensuremath{\mathcalboondox{Y}}}}

\newcommand\U{{\mathcal{U}}} 
\newcommand\Z{{\mathcalboondox{Z}}} 
\newcommand\V{{\mathcalboondox{V}}} 

\newcommand{\topS}{{\ensuremath{\mathcalboondox{S}}}}
\newcommand{\topM}{{\ensuremath{\mathcalboondox{M}}}}
\newcommand{\topG}{{\ensuremath{\mathcalboondox{G}}}}
\newcommand{\topW}{{\ensuremath{\mathcalboondox{W}}}}

\newcommand{\topL}{{\ensuremath{\mathcalboondox{L}}}}
\newcommand{\topLulim}{{\ensuremath{\topL{}_u}}}
\newcommand{\topLlim}{{\ensuremath{\topL{}_l}}}

\newcommand{\topU}{{\ensuremath{\mathcalboondox{U}}}}
\newcommand{\topUp}{{\ensuremath{\topU^{(+)}}}}

\newcommand{\topUpm}{{\ensuremath{\topU^{(+)/-}}}}

\newcommand{\topD}{{\ensuremath{\mathcalboondox{D}}}}
\newcommand{\topDp}{{\ensuremath{\topD^{(+)}}}}
\newcommand{\topDm}{{\ensuremath{\topD^-}}}
\newcommand{\topDpm}{{\ensuremath{\topD^{(+)/-}}}}

\newcommand{\topE}{{\ensuremath{\mathcalboondox{E}}}}
\newcommand{\topEp}{{\ensuremath{\topE^+}}}
\newcommand{\topEm}{{\ensuremath{\topE^-}}}
\newcommand{\topEpm}{{\ensuremath{\topE^{+/-}}}}
\newcommand{\topEmp}{{\ensuremath{\topE^{-/+}}}}

\newcommand{\topY}{{\ensuremath{\Upsilon}}}

\newcommand{\topYpm}{{\ensuremath{\topY^{+/-}}}}

\newcommand{\topT}{{\ensuremath{\Theta}}}

\newcommand\pbcc{\textsc{pb-cc}}
\newcommand\ysf{\textsc{sf-y}}

\newcommand\inv{\ensuremath{^{-1}}}
\newcommand\p{\ensuremath{^\prime}}
\newcommand\R{\mathbb{R}}
\newcommand\Rplus{\mathbb{R}_{\geq 0}}
\newcommand\timeset{\mathbb{T}}

\newcommand{\settrimmingFlow}{%
    \setkeys{Gin}{%
        trim = 2cm 0cm 5cm 2cm, clip=true
    }
    \presetkeys{Gin}{clip}{}
}

\newcommand{\settrimmingRegions}{%
    \setkeys{Gin}{%
        trim = 5.5cm 0 4cm 5cm, clip=true
    }
    \presetkeys{Gin}{clip}{}
}

\binoppenalty=\maxdimen
\relpenalty=\maxdimen

\begin{document}
\title{From lakes and glades to viability algorithms: Automatic classification of system states according to the Topology of Sustainable Management}
\author{Tim Kittel\inst{1} \and
Finn M\"{u}ller-Hansen\inst{1, 2}\fnmsep\thanks{Correspondence to Finn Müller-Hansen, \email{mueller-hansen@mcc-berlin.net}} \and
Rebekka Koch\inst{3} \and
Jobst Heitzig\inst{1} \and
Guillaume Deffuant\inst{4} \and
Jean-Denis Mathias\inst{4} \and
Jürgen Kurths\inst{1,5}}
\institute{Potsdam Institute for Climate Impact Research (PIK),
	Member of the Leibniz Association, FutureLab on Game Theory and Networks of Interacting Agents,
	P.O. Box 60 12 03, D-14412 Potsdam, Germany \and
Mercator Research Institute on Global Commons and Climate Change (MCC),
	EUREF Campus 19, Torgauer Stra{\ss}e 12-15,
	10829 Berlin, Germany \and
	 Institute of Theoretical Physics (ITFA), University of Amsterdam (UvA), Science Park 904, 1098 XH Amsterdam, The Netherlands \and
	Université Clermont Auvergne, INRAE, UR LISC, Aubière, France \and
	Institut für Physik, Humboldt-Universität zu Berlin, 
	Newtonstraße 15,
	12489 Berlin, Germany}

\abstract{
The framework Topology of Sustainable Management by Heitzig et al. (2016) distinguishes qualitatively different regions in state space of dynamical models representing manageable systems with default dynamics. In this paper, we connect the framework to viability theory by defining its main components based on viability kernels and capture basins. This enables us to use the Saint-Pierre algorithm to visualize the shape and calculate the volume of the main partition of the Topology of Sustainable Management. We present an extension of the algorithm to compute implicitly defined capture basins.
To demonstrate the applicability of our approach, we introduce a low-complexity model coupling environmental and socioeconomic dynamics. With this example, we also address two common estimation problems: an unbounded state space and highly varying time scales.
We show that appropriate coordinate transformations can solve these problems.
It is thus demonstrated how algorithmic approaches from viability theory can be used to get a better understanding of the state space of manageable dynamical systems.
} 

\maketitle

\def\jh#1{\textcolor{red}{#1}}
\def\fmh#1{\textcolor{blue}{#1}}

\section{Introduction}

Charting pathways to a sustainable future for humanity needs to account for different biophysical and social constraints. First, the concept of \emph{planetary boundaries} developed by \citet{Rockstrom2009} and extended by \citet{Steffen2015} defines a set of biophysical boundaries that ensure human development under Holocene-like conditions on Earth. They set limits to a \emph{safe operating space} in terms of the planet's capacity to accommodate humanity's metabolism. \citet{Raworth2012} extended the framework to include \emph{social foundations}, i.e. a decent livelihood for every member of society enabling just and equitable human development. Combined with the planetary boundaries, they define the \emph{safe and just operating space}, within which a sustainable development is possible.

There is much research on refining the current definitions of the planetary boundaries, e.g.\ for freshwater \citep{Gerten2013} and phosphorus \citep{Carpenter2011}, on downscaling \citep{Hayha2016} and on extending them, e.g.\ to the terrestrial net primary plant production \citep{Running2012}.
Only a limited number of studies focus on their interaction due to the system's intrinsic dynamics \citep{Anderies2013,Heck2016}.

\citet{Heitzig2016} introduced a mathematical framework called the \emph{Topology of Sustainable Management} (\tsm{}) to study the interactions of boundaries and their consequences for a dynamical system.
The framework can be used to analyze models that distinguish between a default dynamic and one or several management options.
This reflects the observation that social-ecological systems are often governed by relatively stable rules.
To change the resulting dynamics, an extra effort is necessary. Also, the dynamics can only be influenced to a certain degree by active management. Complete control is usually not possible.
The boundaries and the dynamics induce a possibly complex partition of the model's state space into different regions, corresponding to a hierarchy of safety levels.
The regions differ qualitatively in how secure they are and how much management they need to either stay within a desirable region or reach it.
The \tsm{} framework allows for the qualitative analysis of a model and thus precedes other analyses like quantitative optimization.

In this paper, we introduce a variant definition of \tsm{} based on \emph{viability theory}, a subfield of control theory. This allows us to operationalize the \tsm{}-framework for computations. Aubin and his collaborators \citep{Aubin1990,Aubin1991,Aubin2011,Aubin2012} developed viability theory to address the problem of maintaining a dynamical system within a set of desirable states.
Viability theory has been applied in many domains, e.g.\ economics \citep{Aubin1997}, wireless sensor node \citep{Kone2017}, language competition \citep{Chapel2010}, social networks \citep{Mathias2017b} and resilience modeling \citep{Martin2004a,Deffuant2011a,Rouge2013}.

There is a growing literature using viability concepts to analyze sustainability in social-ecological models \citep{Delara2008}. For example, \citet{Baumgartner2009} apply viability concepts to define conditions of strong sustainability under uncertainty.
The review by \citet{Oubraham2018} provides about 80 examples using viability in sustainability studies. In this context, viability theory has been applied for example to fisheries \citep{Cury2005,Chapel2008,Doyen2017,Schuhbauer2016}, water \citep{Richter2003}, forests \citep{Mathias2015}, biodiversity \citep{Doyen2013} and climate change \citep{Mathias2017a}. \citet{Mathias2018, Anderies2019} study the effect of knowledge infrastructure on safe operating spaces in stochastic social-ecological systems.
Some research also focuses on the relationship between viability theory and other mathematical concepts of sustainability, for example the Maximin criterion \citet{Doyen2012}.

In contrast to \tsm{}, viability theory does not distinguish between default dynamics and management options.
Even though \citet{Heitzig2016} already indicated relations between \tsm{} and viability theory, they did not spell out the details.
The analysis of the models in \citep{Heitzig2016} was possible to do without computational analysis as they were all two-dimensional. However, to apply the framework to more complex models, we need to automatize the identification of \tsm{} regions in the state space.
Our definition of \tsm{} differs from the original definition because it is based on the notion of reachability as defined in control theory instead of ``safe reachability'' as defined in \citep{Heitzig2016}.
Also, we require a target set be reached in finite but arbitrary large time, instead of infinite time because it simplifies the computation and is more realistic.
For computationally supported estimations of specific models these differences are usually not relevant.

Our operationalization of the \tsm{}-framework is based on two main concepts from viability theory. 
The first one that we use is the \emph{viability kernel}. It is defined as the largest set of states from where it is possible to keep the system's trajectory within a desirable set and thus avoid the transgression of the set's boundaries.
The second main concept is the \emph{capture basin}. It is defined with respect to a \emph{target set}, i.e.\ a chosen set of states that one wants to reach.
The capture basin is the largest set of states from which the target set can be reached.
This concept is closely related to reachability in control theory.

Our variant definition of the \tsm{} framework allows us to apply computational algorithms from viability theory to estimate different regions of the \tsm{} partition. In particular, we use the Saint-Pierre algorithm \citep{Saint1994} for the estimation of viability kernels. This algorithm is well-established and has been extended with classification methods based on machine learning \citep{deffuant2007}, to viability problems based on reachability \citep{Maidens2013} and optimal control problems \citep{Bokanowski2006,Delara2008,Rouge2013}.
We also develop an extension of the Saint-Pierre Pierre algorithm to find specific \tsm{} regions that can only be defined as implicit target sets.
Furthermore, we propose a nonlinear, local time homogenization solving the problem of vastly differing time scales when estimating viability kernels and capture basins using the Saint-Pierre algorithm. It fully homogenizes the time scale of the dynamics, while keeping the major properties of the system invariant.
With the help of these tools, we can visualize the shape and compute the volume of the \tsm{} partition.

We demonstrate the applicability of the proposed methods by analyzing a three dimensional example model, which focuses on the interaction of stylized environmental and socioeconomic processes similar to other social-ecological models \citep{Wiedermann2015,Barfuss2016,Lade2015,Kellie-Smith2011a,Anderies2013,Brander1988,Nitzbon2017}.
As part of the environment, we model the atmospheric carbon stock and its dynamics. The global economic output and a knowledge stock on renewable energy production, that can induce technological change, describe parts of the socioeconomic system.
In our model, there are two boundaries. The first one is the planetary boundary on climate change \citep{Rockstrom2009,Steffen2015} limiting the amount of atmospheric carbon. The second one is a social foundation for prosperity prescribing a lower threshold on yearly economic output.
We analyze the model with two management options: (i) \emph{low growth} and (ii) climate mitigation by inducing an \emph{energy transition}.

We limit the model to key variables to keep its dimensionality low. Thereby, we ensure that it serves to illustrate the power of our automatic \tsm{} estimation.
Furthermore, we design the model to demonstrate that our methods can deal with common model properties that impede a direct application of computational algorithms such as vastly differing time scales and diverging dynamics. A recent study by \citet{strnad2020} used the model to demonstrate an application of deep reinforcement learning for identifying sustainable management strategies.

The remainder of the article is structured as follows. First, we recall the central notions from viability theory in \Cref{sec:viability-theory}. These are used to derive the variant definition of \tsm{} in \Cref{sec:tsm-variant-definition}. In \Cref{sec:estimation}, we show our means to deal with an unbounded state space and then introduce the nonlinear, local time-homogenization. Next, we shortly recall the Saint-Pierre algorithm and present its extension for implicitly defined target sets. In \Cref{sec:aws-model-description} we introduce the example model with its different managements and analyze it using the tools developed before. We close with a discussion and an outlook for future work.

\section{Viability theory}
\label{sec:viability-theory}

In this section, we shortly introduce the notions of viability theory needed in this article, following \citep{Aubin1991,Aubin2001}.
We start with a time-continuous ($t\in \timeset := \R_+$) dynamical control system 
\begin{align}
\dot{\x} &= f\left(\x, u\right) \label{eq:control-system}
\end{align}
with $\x \in \X = \R^n$, the state space, 
and $u \in \U$, the set of all possible values for the control parameter $u$. We call $f$ the right-hand side function (RHS). Note that no dependency of $\U$ on $\x$ is assumed for ease of notation, but the following could be extended for such a case in a straightforward manner.

A function $q_{\x_0}\colon \timeset \rightarrow \X$ is called a \emph{solution} for an (arbitrary) initial condition $\x_0 \in \X$ if there exists a measurable function of time $\pi\colon \timeset \rightarrow \U$, called \emph{policy}, such that for any time $t \in \timeset$ the condition $\frac{dq}{dt}(t) = f(q(t), \pi(t))$ is fulfilled and $q_{\x_0}(0) = \x_0$.

A \emph{viable set} $\V$ of a  \emph{constraint set} $\Y \subseteq \X$ is then defined as a set of initial conditions $\x$ for which there exists a \emph{viable solution} $q_\x{}$ that stays within $\Y$ forever:
\begin{align}
\forall \x{} \in \V{}\ \exists q_{\x}\ \forall t \in \timeset \colon q_{\x}(t) \in \Y .
\end{align}
The largest viable set of $\Y$ is called the \emph{viability kernel} $\text{Viab}_\U(\Y)$.
The set of possible controls $\U$ is given as a subscript as we will distinguish different controls later.
If one assumes only some constant control $u_f \in \U$, i.e.\ any solution has a policy $\pi(t) = u_f$, the corresponding viability kernel is called a \emph{viability niche}
\begin{align}
\text{VN}_{u_f}(\Y) :=  \text{Viab}_{\{u_f\}}(\Y) .
\end{align}

The \emph{capture basin} of a \emph{target set} $\Z \subseteq \X$ is the subset of the state space in $\X$ for which there exists a solution of \Cref{eq:control-system} that reaches $\Z$ in finite time
\begin{align}
\begin{aligned}[l]
\text{Capt}^\Y_\U \left(\Z\right) =  & \left\lbrace \x_0 \in \X ~|~ \exists \, \text{ solution}\, q_{\x_0}\ \exists T \in \timeset \colon \right. \\
& 
 \left. q_{\x_0}(T) \in \Z \wedge \forall t < T\colon q_{\x_0}(t) \in \Y \right\rbrace .
\end{aligned}  
\end{align}
In case no constraint set is given, the constrained set is assumed to equal the entire state space, i.e.\ $\text{Capt}_\U \left(\Z\right) :=\text{Capt}^\X_\U \left(\Z\right)$.

\section{A variant of the topology of sustainable management based on viability theory}
\label{sec:tsm-variant-definition}

\begin{figure}
\centering
\includegraphics[width=\columnwidth]{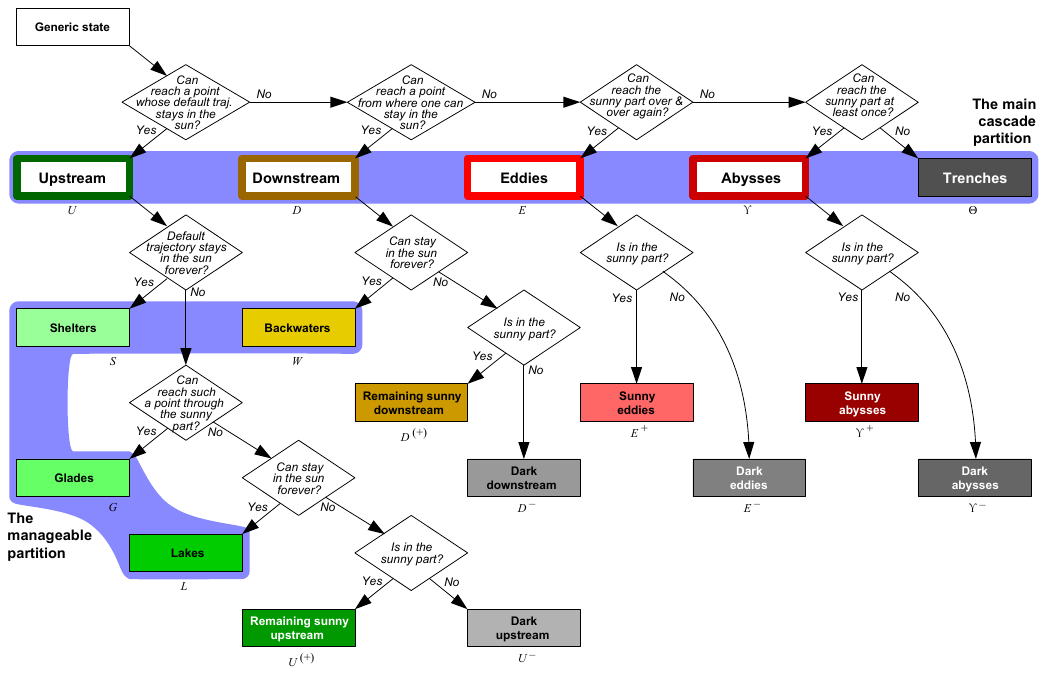}
\caption{The partition of the state space from the Topology of Sustainable Management (reprinted from Heitzig et al. \citep{Heitzig2016}). The decision tree shows to which region a state belongs depending on the need to use the management options and the reachability of the desirable ''sunny`` region (main cascade) or the shelter (manageable partition).}
\label{fig:decision-tree}
\end{figure}

In this section, we present a variant definition of the \tsm{}-partition based on viability theory.
As viability theory, \tsm{} builds on the distinction between a desirable and an undesirable region (set), called ''sunny`` and ''dark`` regions in the \tsm{} framework. However, it additionally distinguishes between default dynamics and management options in the system.
This distinction induces a fine partition of the state space into qualitatively different regions from which the desirable region can or cannot be reached. This also gives rise to multiple dilemmas.
The decision tree in \Cref{fig:decision-tree} shows the different regions of the \tsm{}-partition. A state belongs to a specific region depending on whether it is feasible to stay within it or reach the desirable region and whether management options are needed to do so.
From the perspective of viability theory, \tsm{} can thus be seen as an extension building on the distinction of controls between default and management.
For the full framework, we refer the reader to Ref.~\citep{Heitzig2016}.

\begin{figure}
\centering
\includegraphics[width=\columnwidth]{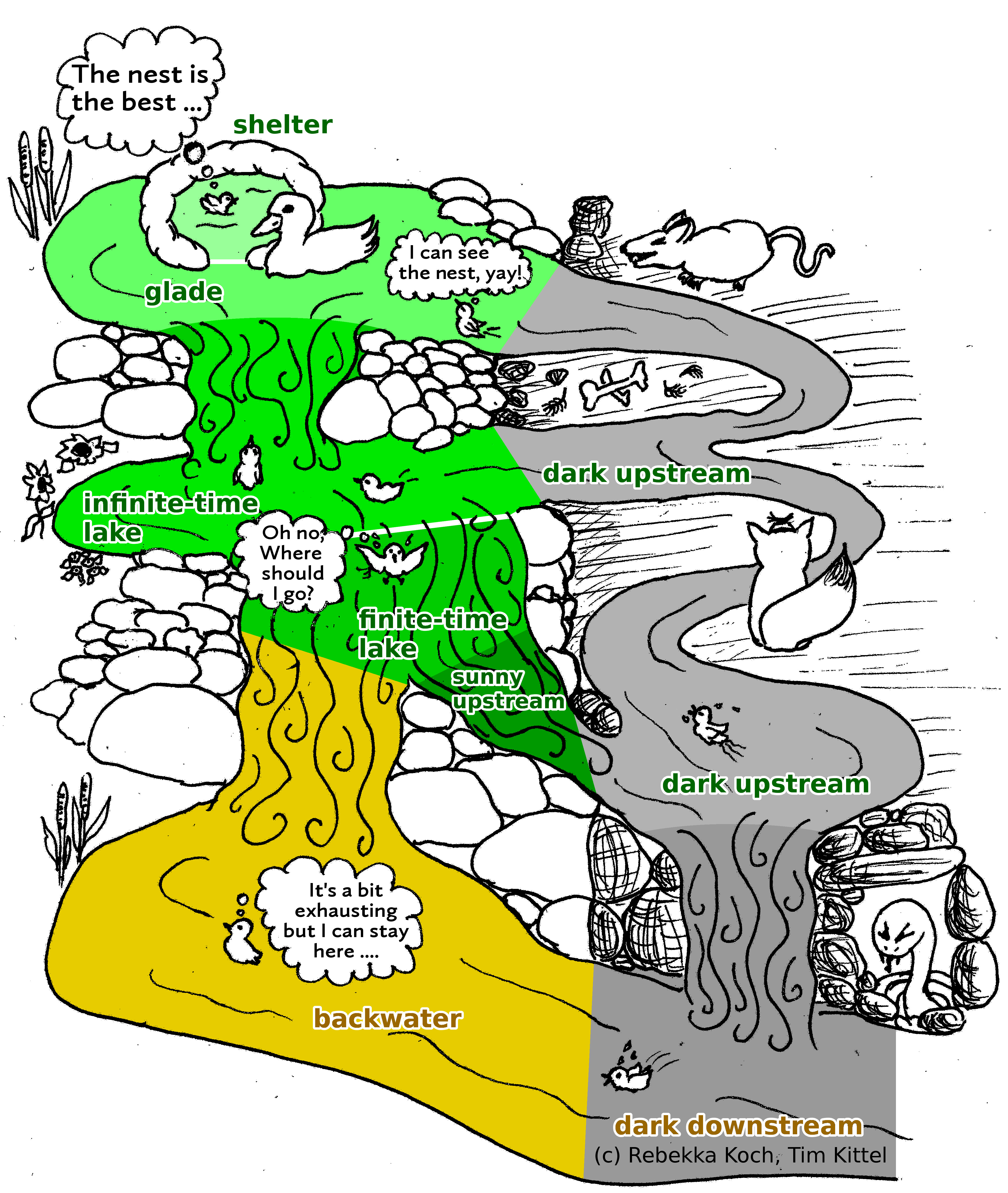}
\caption{Illustration of qualitatively different regions from the \tsm{} framework introduced in \Cref{sec:tsm-variant-definition}. The river flow represents the default dynamics of the system. State space is divided into a desirable (left) and an undesirable region (right). The ducklings, representing different possible states of the system, can either float with flow or actively swim (management option). In the nest, which corresponds to the \emph{shelter}, the ducklings do not have to swim and can stay in the desirable regions without effort. Outside they will slowly drift down or need to swim against the downward-flowing stream to stay in the desirable region. In areas with curly stream lines, the flow is so strong that the ducklings cannot move against it.
}
\label{fig:duckpic}
\end{figure}

We recapture some central notions of \tsm{} using an illustrative example of ducklings in \Cref{fig:duckpic}. The water region represents the state space and the streamlines represent the dynamics. The ducklings can either swim with the flow (default dynamics) or struggle and swim against it (management). However, it is not possible to swim up a waterfall once they have dropped down. The desirable region, a safe environment providing enough food and nesting places for the ducklings, is on the left and the undesirable region full of predators on the right.   
In the following, we introduce the general concept for each region and then illustrate it with the ducklings example.

Our definitions of \tsm{}-regions are based on a general control system as \Cref{eq:control-system}, for which we additionally require that a \emph{default control} $u_0 \in \U$ is separated out from the set of \emph{all possible controls} $\U$. $\U_m = \U - \{u_0\}$ is the set of \emph{manageable controls}.\footnote{Here and in the following we use the lax difference and union notation with ``$-$'' and ``$+$'' for sets.}
Hence we call $f(x, u_0)$ the \emph{default flow} giving rise to the \emph{default dynamics} and the dynamics corresponding to the manageable controls is called \emph{management options}.
Furthermore, we require a division of the state space into \emph{desirable} $\Xp \subseteq \X$ and \emph{undesirable} $\Xm := \X-\Xp$ regions.

As already sketched in \citep{Heitzig2016}, the viability perspective allows us to base the partition of state space upon two basic regions, the \emph{shelters} and the larger \emph{manageable region}.
In the shelters $\topS$, the system will stay in the desirable region forever without management, i.e. only by following the default dynamics described by the default control $u_0$. Thus, they are the safest regions in state space. In terms of viability theory, they are the viability niche of $\Xp$ under $u_0$,
\begin{align}
&\topS:= \text{VN}_{u_0} \left( \Xp \right).
\end{align}
For the ducklings in \Cref{fig:duckpic} the shelter is the nest, a place in which they can stay in a safe environment forever without swimming against the stream. 

The manageable region $\topM$ is the part of the desirable set where one can stay by making use of all management options. With that, it is the viability kernel of $\Xp$ when taking into consideration the full set of controls $\U$.
\begin{align}
&\topM:= \text{Viab}_\U \left( \Xp \right).
\end{align}
In \Cref{fig:duckpic}, this corresponds to all parts of the river in the desirable region on the left. However, there can be other parts of the desirable region that do not allow to stay in the desirable region forever, even with management (``remaining sunny downstream'', not depicted in \Cref{fig:duckpic}).

In contrast to the definition here, \citet{Heitzig2016} define the shelters as the ``invariant open kernel'' of $\Xp$, which is slightly smaller than $\text{VN}_{u_0} \left( \Xp \right)$ because it does not contain its boundary. Likewise, they define the manageable region in terms of the ``sustainability kernel'' instead of the viability kernel, which again is a slight subset that does not contain the boundary. Analogous remarks apply to all the other regions defined below.

In \tsm{}, the \emph{upstream} is defined as the region from where it is possible to reach the shelter(s) with or without management. This can include leaving the desirable region for a limited time period. In viability theory, the capture basin formalizes the idea of being able to reach a certain set. Hence, the upstream is the capture basin of the shelter,
\begin{align}
\topU:=\text{Capt}_\U \left( \topS \right)
\end{align}

The \emph{downstream} are those regions from where the manageable region is reachable but it is not possible to reach the shelter(s). From the downstream it is possible to reach a region where one can stay within the desirable region forever using management. The downstream region is therefore the capture basin of the manageable region without the upstream,
\begin{align}
\topD := \text{Capt}_\U \left( \topM \right) - \topU .
\end{align}
\tsm{} partitions the state space into non-intersecting regions. This is why our definitions often need to remove certain sets so they are not classified twice. In the above definition, the removal of the upstream ensures that the shelter cannot be reached from the downstream.

As both the upstream and the downstream may have qualitatively different dynamics inside, \tsm{} introduces a finer subdivision of these regions.
The upstream is further divided into glades, lakes and the remaining sunny and dark upstream.
\emph{Glades} are the region from where it is possible to reach the shelter without leaving the desirable set. This means there is only a temporary need for management to stay in the desirable region forever. After reaching the shelter, the system will remain in the desirable set without management. In terms of viability theory, the glades are the capture basin of the shelter (excluding the shelter itself) when only considering the sunny region,
\begin{align}
\topG = \text{Capt}_\U^\Xp \left( \topS \right) - \topS .
\end{align}
In our ducklings example, the water around the nest is a glade.

From the \emph{lakes}, the shelters are only reachable trough the undesirable region. Thus, they are the intersection between the upstream (shelter is reachable) and the manageable region with the shelters and the glades removed,
\begin{align}
\topL = \topU \cap \topM - \topS - \topG. 
\end{align}  
Being located in a lake imposes a qualitative choice, the \emph{lake dilemma}: Inside a lake one can choose between eventual safety and uninterrupted desirability. Deciding for the first requires crossing the undesired region to finally reach the shelter. But deciding for the second implies that there will always be a need for management options.

As a further refinement, let us distinguish two kinds of lakes, \emph{time-limited lakes} \topLlim{} and \emph{time-unlimited lakes} \topLulim{}, which can be defined as
\begin{align}
\topLulim{} &= \text{Viab}_\U \left(\topL\right),\\
\topLlim{} &= \topL - \topLulim.
\end{align}
Because time-unlimited lakes are the viability kernel of the lakes, there is no time pressure for a decision on the dilemma.
In contrast, there is a fixed deadline in time-limited lakes. This leads to a stronger form of the dilemma called the \emph{pressing lake dilemma}. In \Cref{fig:duckpic}, the time-unlimited lake is the waterfall starting from the glade and the subsequent calmer waters. The ducklings can swim in the calmer waters as long as they want to.
The time-limited lake is the following waterfall that splits into two streams. Here, the ducklings drop down the waterfall and they have only a moment to decide for the left or the right.

The rest of the upstream is split into the \emph{remaining sunny} $\topUp$ and the \emph{dark upstream} $\topUp$, depending on whether a state is in the desirable region or not. They can be defined as the intersection of the upstream without the manageable region with the desirable or undesirable region,
\begin{align}
\topUpm = (\topU - \topM) \cap \Xpm.
\end{align}

The downstream is subdivided into backwaters, remaining sunny and dark downstream. The \emph{backwaters} are those regions where it is possible to remain in the desirable region using management but the shelters are not reachable. The first condition means that they are a part of the manageable region and the second condition excludes that it is in the upstream, thus
\begin{align}
\topW = \topM - \topU.
\end{align}
The backwaters can also be defined as the intersection of the manageable region and the downstream, $\topW = \topM \cap \topD$. The calm waters in the left lower corner of \Cref{fig:duckpic} is a backwater because the ducklings can swim against the stream and stay inside the backwater but due to the waterfalls, they cannot reach the nest.
 
Within the capture basin of the backwaters, there are no further regions such as the lakes in the upstream. Analogously to the upstream, the rest of the downstream is divided into the \emph{remaining sunny} $\topDp$ and the \emph{dark downstream} $\topDm$, depending on whether they are in the desirable region or not. These regions are therefore the intersection of the downstream without the backwater with the desirable or undesirable region, respectively,
\begin{align}
\topDpm = \left( \topD - \topW \right) \cap \Xpm . 
\end{align}

If the desirable region can be reached over and over again one is inside the \emph{eddies} $\topE$ that are divided into \emph{sunny eddies} $\topEp$ and \emph{dark eddies} $\topEm$ (these and the following regions are not depicted in \Cref{fig:duckpic}).
The metaphorical image behind the naming is that of a circular flow of which one segment is in the desirable and the other one is in the undesirable region. 
They are the maximal pair of sets fulfilling
\begin{subequations}
\begin{align}
	&\topEpm \subseteq \Xpm - \topU - \topD , \label{eq:eddies-1} \\
	&\topEpm \subseteq \text{Capt}_\U (\topEmp) , \label{eq:eddies-2} \\
	&\topE = \topEp + \topEm . \label{eq:eddies-3}
\end{align}
\end{subequations}
\Cref{sec:estimation-eddies} and \Cref{app:existence-eddies} provide details why eddies can be defined like this.
\emph{Trenches} are regions from which the desirable region cannot be reached ever again. In terms of viability theory, they are the regions that do not belong to the capture basin of the desirable region,
\begin{align}
\Theta = \X - \text{Capt}_\U \left( \Xp \right).
\end{align}

To complete the partition of the state space, the \emph{abysses} are defined as the rest of the state space not belonging to any region mentioned before,
\begin{align}
\topY= \X - \topU - \topD - \topE - \Theta.
\end{align}
Inside the abysses, on can reach the desirable region a finite number of times only. There are sunny and dark abysses $\topYpm = \topY \cap \Xpm$.
This completes the main partition of the \tsm{} framework.

\section{Estimation}
\label{sec:estimation}

To estimate the \tsm{}-partition, we use the definitions derived in \Cref{sec:tsm-variant-definition} and the established Saint-Pierre algorithm \citep{Saint-Pierre2001} to find the viability kernels and capture basins that our definitions build upon.
We apply the Saint-Pierre algorithm because it is the standard algorithm in the field and operationalizes the most fundamental idea how to estimate viability kernels.
Furthermore, we chose it as a first step for automatic TSM region detection because it is comparatively simple to implement.
The Saint-Pierre algorithm always converges to the viability kernel under well-defined conditions and when decreasing adequately the time step and the space resolution.
Other algorithms may be more efficient on specific examples but the conditions of their convergence to the viability kernel are more restrictive \citep{Krawczyck2011, Bokanowski2006, deffuant2007, Coquelin2007}.

The Saint-Pierre algorithm is based on finitely discretizing the state space and then using local linear approximations of the dynamics. Hence, it is applicable to bounded state spaces only.
In case the relevant domain of the state space is unbounded, we need to map it to a bounded space first. Also, vastly differing time scales might be problematic for the linear approximations, so there is a need to homogenize the time scales.
In this section, we explain the basic ideas of the Saint-Pierre algorithm and introduce the two coordinate transformations to make it applicable. Finally, we extend the algorithm to estimate implicitly defined capture basins for identifying the eddies of the \tsm{}-partition.

\subsection{The Saint-Pierre algorithm}
\label{sec:saint-pierre-algorithm}

The Saint-Pierre algorithm \citep{Saint1994} was developed in order to estimate the viability kernel of a control system \Cref{eq:control-system}.

It starts with a discretization $\Y_h$ of the constraint set $\Y$ where a point $x \in \Y$ is at most at a distance $h$ of a point $\y \in \Y_h$\footnote{Often, a regular grid with resolution $h$ is chosen for this discretization, but this is by no means a necessity. We use a regular grid because we already apply a coordinate transformation (see \Cref{sec:unbounded-state-space-general}) that best resolves the scales of interest in our model.}. Furthermore, a small time step $\Delta t > 0$ is chosen and it supposes that the set of controls $\U$ is discrete (if not, it is also discretized). It supposes that $f$ is $l$-Lipschitz and there exists an upper bound $M$ of $\Y$. 

For each point $\x \in \Y_h$ and for each control $u \in \U$, the algorithm computes the successors $S(\x, u)$ of $x$ when applying control $u$ for a linearized, extended dynamics defined from $f$. Note that the linearization is done in the $\x$ variable only. It is extended in the sense that the successors include all the points located in a ball around $\x + f(\x,u)\cdot \Delta t$. The successors $S(\x, u)$ of $x$ when applying control $u$ are given by
\begin{equation}
 S(\x,u) = \left\{ \y \in \Y_h ~|~ \left\|\y - (\x +  f(\x,u)\cdot \Delta t)\right\| \leq h + \frac{Ml}{2}\cdot \left(\Delta t\right)^2 \right\}.
\end{equation}
This extension of the dynamics guarantees that the algorithm described below converges to the actual viability kernel when $\Delta t$ and the resolution of the grid decrease to 0. Computing and storing all the successors for each point of the grid rapidly becomes computationally heavy when the dimensionality of the the state space is large and the grid resolution is small (this is an example of the famous curse of dimensionality).

Then, the algorithm builds a series of discrete sets (subsets of the grid $\Y_h$) $K_0 = \Y_h, K_1,..., K_n$ such that $K_{i+1} \subset K_i$, defined as follows:

\begin{equation}
K_{i+1} = \left\{ x \in K_i ~|~ \exists u \in \U \colon S(x,u) \cap K_i \neq \emptyset  \right\}
\end{equation}

After a finite number of steps, the algorithm reaches a fixed point, i.e.\ $K_{n+1} = K_n$. The set $K_n$ is the viability kernel of $\Y_h$ for the linearized, extended discrete dynamics. \citet{Saint1994} shows that this set converges to the viability kernel of the continuous time dynamics when $\Delta t$ and $h$ tend to 0 appropriately. Note that the approximations are done from the exterior of the viability kernel: generally, the approximation includes points that do not belong to the actual viability kernel (their proportion decreases when the resolution of the grid decreases).

This algorithm has been extended by \citet{deffuant2007} for using continuous sets $K_i$, using a machine learning  algorithm that takes as input the points of the grid that belong to $K_i$ and the ones that do not, and derives an approximation of its boundary. This opens the possibility to represent continuous viability kernels that are defined more conveniently than a huge set of points.

A slight modification of the algorithm described above enables us to approximate the capture basin (see \cite{Saint-Pierre2001} for details).
We start with a discretization of the state space $\X_h$, analogously to $\Y_h$ above, and define the discretized target set $\Z_h = \Z \cap \X_h$ for a target set $\Z$. Again, we create a series of discrete sets $K\p_i$ with $K\p_0 = \Z_h$ and where the successors of all elements in $K\p_{i+1}$ are in $K\p_i$
\begin{equation}
K\p_{i+1} = \left\{ x \in \X_h ~|~ \exists u \in \U \colon S(x,u) \cap K\p_i \neq \emptyset  \right\}.
\end{equation}
Again, after a finite number of steps, the algorithm reaches a fixed point, i.e.\ $K\p_{n+1} = K\p_n$, and $K\p_n$ is the capture basin of $\Z_h$ in $\X_h$ for the linearized extended discrete dynamics. In contrast to the viability estimation, this is an approximation from the interior.

Improvements and extensions to this algorithm are currently under intensive research. Relations to dynamical programming \citep{Frankowska1989} and other extensions \citep{Bayen2002,Saint-Pierre2001} can provide the minimal time to reach the target set. Also, one can even find controllers that drive the system to the target set \citep{Cardaliaguet1998,Lhommeau2007,Chapel2011}.

\subsection{Dealing with an unbounded state space}
\label{sec:unbounded-state-space-general}

There are multiple ways to map an unbounded state space to a bounded one, depending on the specific need for the system. In case of the example system analyzed later, each coordinate is bounded from below and unbounded from above. This is rather common in socioeconomic models, in particular due to unbounded economic growth. Hence, we propose a solution that maps each coordinate separately.

We assume a general dynamical system given by a set of ordinary differential equations
\begin{align}
\dot \x = f(\x)
\end{align}
with $\x \in \Rplus^n$.
In contrast to \Cref{eq:control-system}, there is no control parameter, because the dependence on the control parameter is irrelevant here and the mapping can be done for general ordinary differential equations. We imply an extension to control systems by considering fixed controls.
Then we propose the coordinate transformation 
\bgroup
\renewcommand*{\arraystretch}{1.5}
\begin{align}
\begin{array}{@{}r@{\ }c@{}c@{\enspace}l@{}}
\Phi: & \Rplus^n &\longrightarrow& \left[0, 1\right)^n \\
&(\x_i) &\longmapsto& (\frac{\x_i}{\x_{i,mid} + \x_i})  ,
\end{array}
\end{align}
where $\x_{i,mid} \in \Rplus$ are parameters.
\egroup
Applying this transformation on the dynamics leads to a new set of ordinary differential equations
\begin{align}
\dot \y = F(\y) :=& ((D \Phi \cdot f) \circ \Phi\inv)(\y) \\
=& \frac{(1-\y_i)^2}{\x_{i,mid}} f\left( \frac{\y_i \x_{i,mid}}{1-\y_i} \right) \label{eq:general-coordinate-transformed-system}
\end{align}
for $\y \in [0,1)^n$, i.e.\ inside a bounded space, where $D\Phi$ is the Jacobian of $\Phi$ and $\circ$ is the symbol for function composition.
The parameters $x_{i,mid}$, summarized to the vector $x_{mid}$, define the scales for each coordinate that is resolved best when discretizing the state space for numerical estimations because $\Phi(x_{mid}) = (\frac{1}{2}, \frac{1}{2}, \dots)^T$. So they should be taken to be around the main region of interest.

For the dynamical system including the management options, we assume the control parameters to be fixed and let the system switch between the default dynamics and the management options instantaneously. Therefore, the transformation given in this section can also be applied to the controlled system.

\subsection{Nonlinear, local time-homogenization}
\label{sec:time-homogenization}

A problem during the estimation of viability kernels from \Cref{sec:viability-theory} is a possibly inhomogeneous time scale, i.e.\ that the (norm of the) RHS function of the control system \eqref{eq:control-system} can have values through \emph{several orders of magnitude}. For instance, models like \Cref{eq:general-coordinate-transformed-system} often lead to divergences at the upper boundary of a coordinate.

This problem can be addressed by rescaling the time of the system in a nonlinear way. The used definitions of viability theory depend only qualitatively but not quantitatively on time: It does not matter for the definition of viability kernels and capture basins how long the system takes to reach the target set or desirable region, it only matters that it is possible to reach in a finite time period. Hence, viability kernels and capture basins in the rescaled system are equivalent to the original ones. Because there is no functional difference between these quantities in the original and the homogenized system, it is technically possible to use the latter in order to estimate them in the former.

As the control parameter is not necessary for the rescaling of the system, we use a differential equation
\begin{align}
\dot \y = F(\y) \label{eq:original-system-for-time-homogenization}.
\end{align}

We propose to use the new system
\begin{align}
	\dot \y = \tilde F(\y) := \frac{F(\y)}{\|F(\y)\| + \epsilon}. \label{eq:time-rescaled-system}
\end{align}

Assuming $\epsilon$ is small enough, this new system generally fulfills three criteria:
\begin{enumerate}
	\item The systems \eqref{eq:original-system-for-time-homogenization} and \eqref{eq:time-rescaled-system} are \emph{orbitally equivalent} (cf.\ \cite[p. 42]{Kuznetsov1998}, Definition 2.4), i.e.\ the trajectories of solutions with the same initial conditions follow the same path. In other words, only the time has been rescaled.
	\item Everywhere away from fixed points $\| F(\y) \| \gg \epsilon$ holds and hence the time scale is properly homogenized 
	\begin{align}
	\| \tilde F(\y) \| = \frac{ \| F(\y) \| }{ \| F(\y) \| + \epsilon} \approx  1	.
	\end{align}
	
	\item Within a small enough environment of a fixed point, $ \| F(\y) \| \ll \epsilon$ holds. We can therefore approximate
	\begin{align}
	\tilde F(\y) = \frac{ F(\y) }{ \| F(\y) \| + \epsilon} \approx \frac{ 1}{ \epsilon} F(\y).
	\end{align}
	This implies that the function $\tilde F$ goes to zero at fixed points of the original system with the same properties as $f$ at these points (e.g. $\mu$-Lipschitz or $C^\infty$, same Lyapunov-Exponents etc.).
\end{enumerate}

Because the units of the coordinates of $\y$ might be different from each other, there is no real physical interpretation of $\tilde F$. But that is not necessary either as it is only an auxiliary system for the estimation with the Saint-Pierre Algorithm.

\subsection{Estimation of implicitly defined capture basins: eddies}
\label{sec:estimation-eddies}

For capture basins we only care about entering a target set at least once. However, eddies are defined by being able to reach the desirable region again and again, so the most natural definition is an implicit one as in \Cref{eq:eddies-1,eq:eddies-2}. In order to estimate them, we find an alternate definition in terms of a limit process.

We start by defining the largest sets that could contain eddies
\begin{subequations}
\begin{align}
\topEp_0 &= \Xp - \topU - \topD, \\
\topEm_0 &= \Xm - \topU - \topD,
\end{align}
\end{subequations}
and then use the iteration step
\begin{subequations}
\begin{align}
	\topEm_i &= \text{Capt}_\U (\topEp_{i-1}) \cap  \topEm_{i-1}, \label{eq:eddies-iteration-step-minus}  \\
	\topEp_i &= \text{Capt}_\U (\topEm_i) \cap \topEp_{i-1} \label{eq:eddies-iteration-step-plus}
\end{align}
\end{subequations}
for $i = 1,2,\dots$. Note that $\topEm_i$ used in \Cref{eq:eddies-iteration-step-plus} is already computed in \Cref{eq:eddies-iteration-step-minus}. So this can really be seen as a step-by-step prescription.
Thus, the eddies can be recovered as 
\begin{subequations}
\begin{align}
	\topEp &= \lim_{i\to\infty} \topEp_i, \label{eq:eddies-iteration-limit-plus} \\
	\topEm &= \lim_{i\to\infty} \topEm_i. \label{eq:eddies-iteration-limit-minus}
\end{align}
\end{subequations}
The limit exists because both sequences are monotone and nonincreasing.
The Saint-Pierre algorithm works on a discretized state space with finitely many elements. Hence, there exists an $k\in \mathbb{N}$ such that $\topEp = \topEp_k = \topEp_{k-1}$ and $\topEm = \topEm_k = \topEm_{k-1}$. Thus, the existence of the limit results in the convergence of the algorithm after a finite number of steps.

This iteration process follows the idea of being able to visit the desirable region again and again and is an algorithmic description for the estimation of eddies. Other similarly implicitly given sets can be estimated by adjusting this basic idea.


\section{Example: the AYS low-complexity model of climate change, wealth and energy transition}
\label{sec:aws-model-description}

\begin{figure}
\centering
\includegraphics[width=0.8\columnwidth]{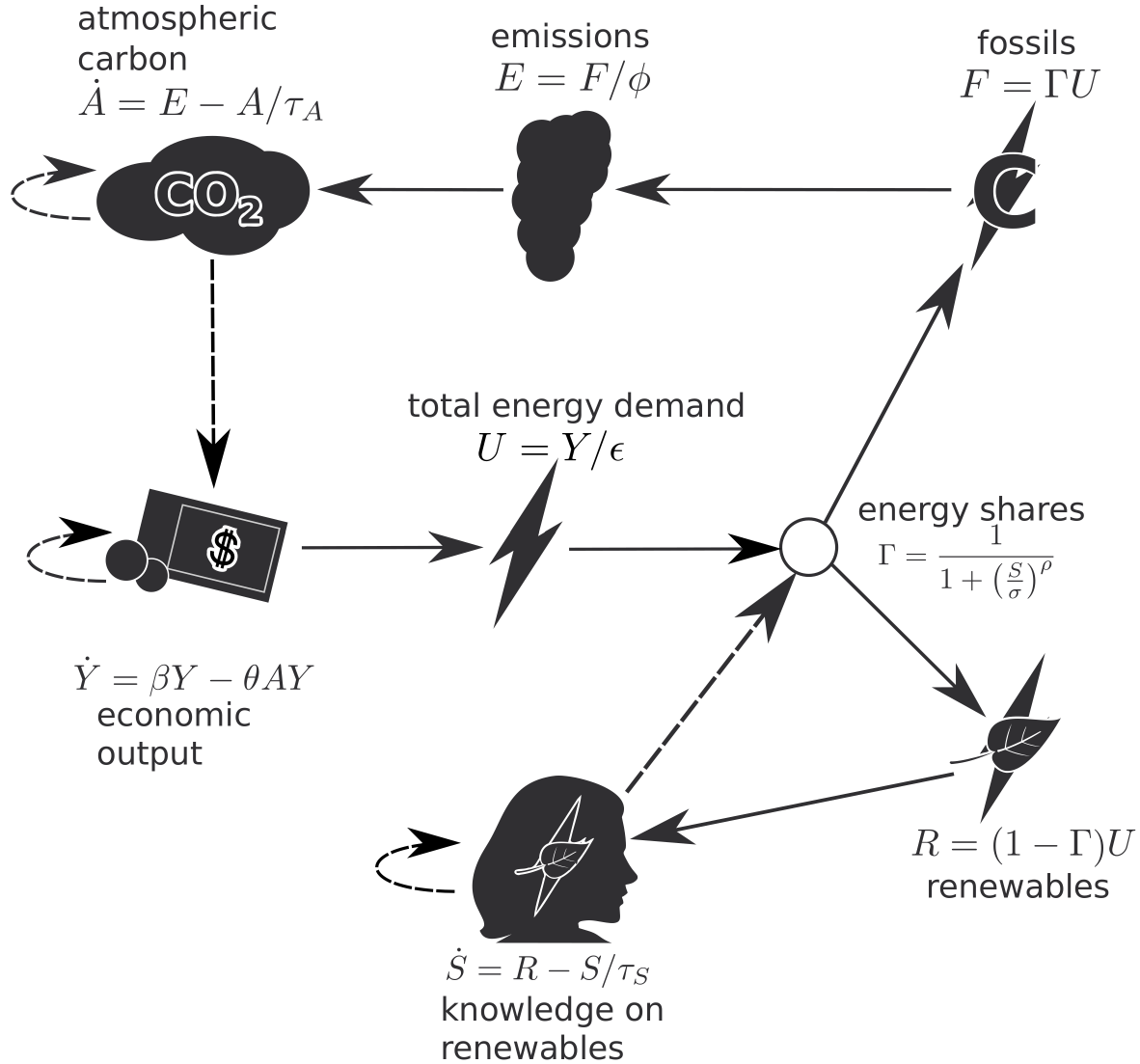}
\caption{The interplay of the three dynamical variables excess atmospheric carbon stock $A$, economic production $Y$ and renewable energy knowledge stock $S$ and the five dependent variables energy demand $U$,
fossil energy flow $F$, renewable energy flow $R$,
emissions $E$ and share of the fossil sector $\Gamma$. }
\label{fig:aws-model-graph}
\end{figure}

Here, we demonstrate the applicability of our operationalization of the \tsm{}-framework using a three dimensional example model. We start with describing the model's default dynamics, the management options and the desirable region. Then, we analyze the model's attractors and estimate the current state within the state space. Finally, we apply the coordinate transformations and discuss the resulting \tsm{} regions found with the Saint-Pierre algorithm.

\subsection{Model default dynamics}

To develop our low-complexity example model, we took inspiration from \citet{Kellie-Smith2011a}
and added a renewable energy sector with a learning-by-doing dynamics.
The model describes the interactions between the atmosphere and a stylized economy whose energy supply comes from either fossil fuels and renewable resources. While the atmosphere reacts to CO$_2$ emissions by increasing temperatures, the higher temperatures have negative impact on economic output via climate damages. The share between the two energy sources is determined by their relative cost, which depends on previous use of renewables and optional policies.
The model structure is depicted in \Cref{fig:aws-model-graph}. It comprises three state variables.

1) The first variable is the \emph{excess atmospheric carbon stock} $A$ [GtC = giga tons of carbon],
measured w.r.t.\ a pre-industrial level $A_0 \approx 600\,$GtC. 
It increases with current CO$_2$ emissions $E$ [GtC/a = GtC per year]. 
Taking $A_0$ as an estimate for the long-term no-emissions equilibrium value, 
we assume $A$ approaches zero if $E = 0$, due to carbon uptake by oceans, plants and soil. 
To keep the complexity of the model as low as possible, 
we do not explicitly model a carbon cycle as in \citet{Anderies2013}
but simply assume the carbon uptake leads to an exponential relaxation towards equilibrium 
on a characteristic time scale of $\tau_A \approx 50\,$a [a = years].
Hence our first model equation is
\begin{align}
	\frac{dA}{dt} &= E - A / \tau_A, \label{eq:derivate-A}
\end{align}
where $E$ will be derived below from economic assumptions.

2) The second variable is \emph{economic output} or production $Y$ [US\$/a] of the global economy. It thus represents the gross world product as a common indicator. 
We assume the economy to have a positive basic growth rate $\beta \approx 3\,\%$ [1/a] 
and additional climate impacts as in \citep{Kellie-Smith2011a}. 
As a proxy for temperature we simply use $A$, effectively assuming an infinitely fast greenhouse effect.
Hence this terms is represented by $-\theta A Y$ where $\theta \approx 8.57\cdot 10^{-5}\, /$(GtC\,a) is a temperature sensitivity parameter 
chosen such that the total growth rate $\beta - \theta A$ becomes negative 
when $A$ exceeds the level corresponding to a global warming of $+2\,^\circ\text{C}$.
This gives
\begin{align}
	\frac{dY}{dt} &= \beta Y - \theta A Y. \label{eq:derivative-Y}
\end{align}

3) The third state variable is the \emph{renewable energy knowledge stock} $S$ 
that indicates how much knowledge is available for the production of renewable energy $R$ [GJ/a = giga joule per year].
In accordance with Wright's law (e.g., \citep{Nagy2013}) of learning-by-doing, 
we basically identify $S$ with the past cumulative production of renewables 
and thus measure it in units of [GJ]. 
To account for the human capital component, 
we additionally assume that knowledge depreciates on a characteristic time scale of $\tau_S \approx 50\,$a.
We choose this value to be in the upper range of a typical professional life span.
Cumulation and depreciation then give
\begin{align}
	\frac{dS}{dt} &= R - S / \tau_S,
	\label{eq:derivative-S}
\end{align}
where $R$ will be derived below.

Finally, to determine $E$ and $R$, we use the following simplistic economic assumptions.
The energy demand $U$ [GJ/a] is proportional to the economic output
\begin{align}
U = Y / \epsilon ,
\label{eq:definition-U}
\end{align}
where $\epsilon \approx 147\,$US\$/GJ is an energy efficiency parameter.
This demand is satisfied by a mix of fossil and renewable energy
which are assumed to be perfect substitutes 
(and ignoring other energy sources such as bioenergy).
Their respective shares are determined by a price equilibrium.
We assume convex monomial cost functions
and unit costs of renewable energy that show a power-law decay with growing $S$ \citep{Nagy2013}.
This implies that the fossil sector has a share given by the sigmoidal function 
\begin{align}
	\Gamma &= \frac{1}{1 + \left(\frac{S}{\sigma}\right)^\rho}, \label{eq:definition-gamma}
\end{align}
where $\sigma\approx 4\cdot 10^{12}\,GJ$ is the break-even knowledge level at which renewable and fossil extraction costs become equal, 
and $\rho\approx 2$ is a dimensionless parameter determined from the cost convexity and learning rate. 
$\Gamma$ approaches unity (no renewables) as $S\to 0$ and zero (no fossils) as $S\to\infty$. 
Fossil and renewable energy flows and emissions are then
\begin{align}
	F &= \Gamma\, U, & R &= \left(1 - \Gamma\right) U, & E &= F / \phi, \label{eq:definition-F-R-E}
\end{align}
where our final parameter $\phi \approx 4.7\cdot10^{10}\,$GJ/GtC is the fossil fuel combustion efficiency.
This completes the equations of the AYS model.

\begin{figure}
\centering
\settrimmingFlow
\includegraphics[width=\columnwidth]{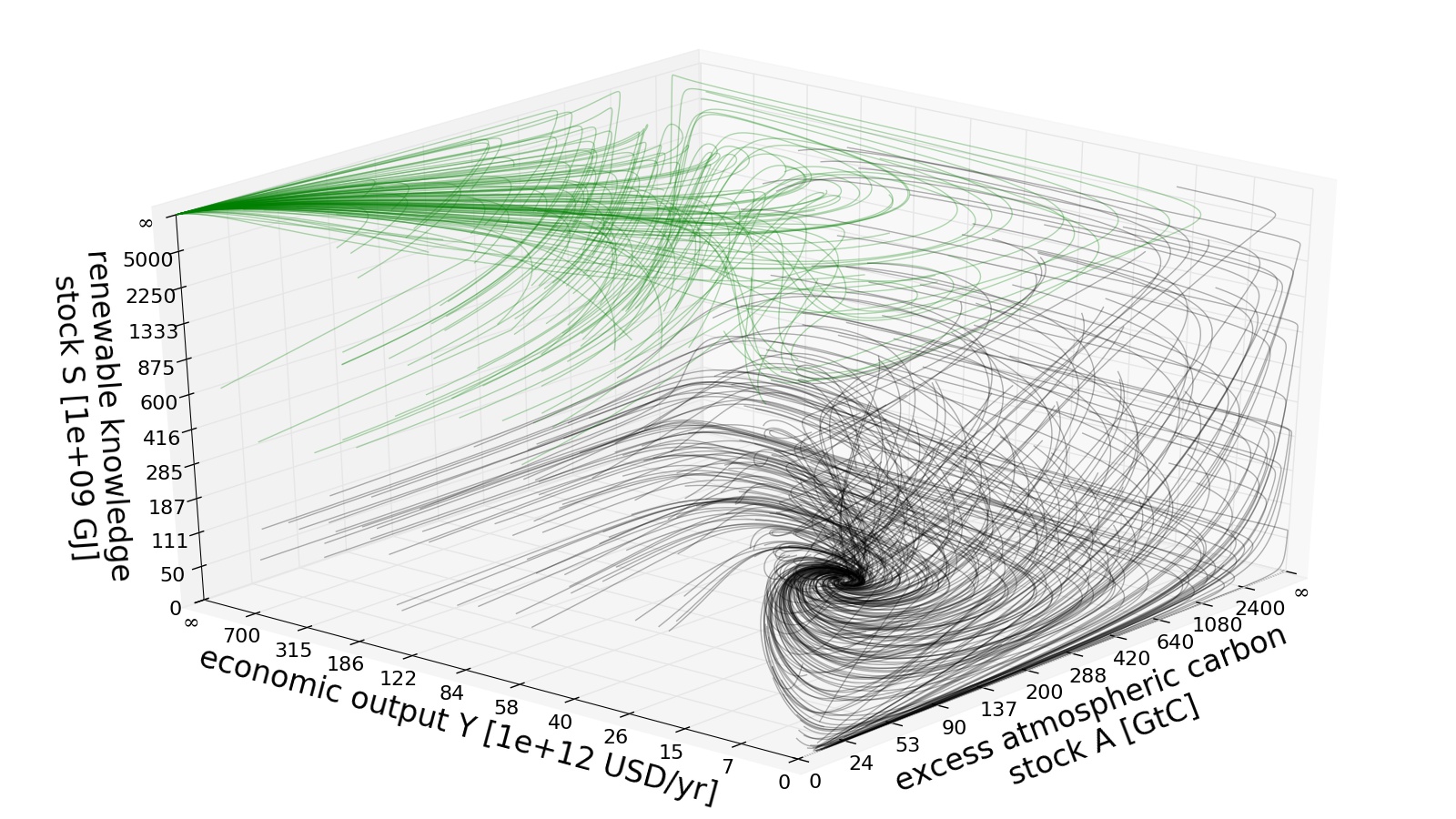}
\caption{The default flow of the AYS model (\Cref{eq:derivate-A,eq:definition-gamma,eq:definition-F-R-E,eq:derivative-Y,eq:derivative-S,eq:definition-U}) is sampled with trajectories from randomly distributed initial conditions on nonlinearly scaled axes so the full states space $X=\Rplus^3$ is displayed (more details in \Cref{sec:aws-compatification}). Green trajectories end up at the green attractor $x_g$ and black ones at $x_b$ (see \Cref{sec:aws-attractors}).}
\label{fig:aws-default-flow-without-boundaries}
\end{figure}

The three dynamical variables $A$, $Y$ and $S$ are interrelated due to the various connecting equations and nonlinearities arise particularly due to \Cref{eq:definition-gamma} and \Cref{eq:derivative-Y}.
The model's default dynamics is also specified by the eight parameters $\tau_A$, $\beta$, $\theta$, $\tau_S$, $\epsilon$, $\rho$, $\sigma$ and $\phi$.
\Cref{app:paramter-estimation} contains details on how we estimated these parameters.
The resulting flow is depicted in \Cref{fig:aws-default-flow-without-boundaries}, where the basins of the system's two attractors (discussed in \Cref{sec:aws-attractors}) are already colored differently.

\subsection{Management options (controls)}
\label{sec:aws-management-options}

The above parameter values define what we consider the {\em default dynamics} because they represent a business-as-usual case. This means humanity applies no specific management that would alter the way things usually go.

In addition to the default dynamics, we consider two distinct management options. They represent possible policy choices that may be combined in any way, leading to more or less shifted trajectories. Because the management options are defined as changes in the parameters $\beta$ and $\sigma$, the control parameter $u$ of the model is two-dimensional and the particular control $u_0$ for the default flow is
\begin{align}
u_0 = \begin{pmatrix}
\beta \\
\sigma
\end{pmatrix}.
\end{align}

1) The first distinct management option we consider is {\em low growth} (\lg{}). This corresponds to policies that voluntarily reduce the basic growth rate $\beta$ to half its value $\beta_\lg = 1.5\%/$a, leading to the control
\begin{align}
u_\lg = \begin{pmatrix}
\beta_\lg \\
\sigma
\end{pmatrix}.
\end{align}

2) The second management option we introduce aims at mitigating climate change by inducing an {\em energy transition (\et{})}, e.g.\ via taxing fossils and/or subsidizing renewable energy use. 
These policy instruments shift the relative costs of fossil and renewable energy, which according to \Cref{eq:definition-gamma} can be effected in our model by a reduction of $\sigma$. Hence, we represent this option by reducing $\sigma$ to approx. $(1/2)^{1/\rho} = 1/\sqrt{2}$ of its default value, i.e.\ to $\sigma_\et = 2.83\cdot 10^{12}$\,GJ, corresponding to dividing the renewable to fossil cost ratio by half and leading to the control
\begin{align}
u_\et = \begin{pmatrix}
\beta \\
\sigma_\et
\end{pmatrix}.
\end{align}
With these two management options, the set of possible managements is given by $\U_m = \left\lbrace u_\lg, u_\et \right\rbrace$.

\begin{figure}
	\centering
	\settrimmingFlow
	\subfloat[Default dynamics]{\includegraphics[width=0.8\columnwidth]{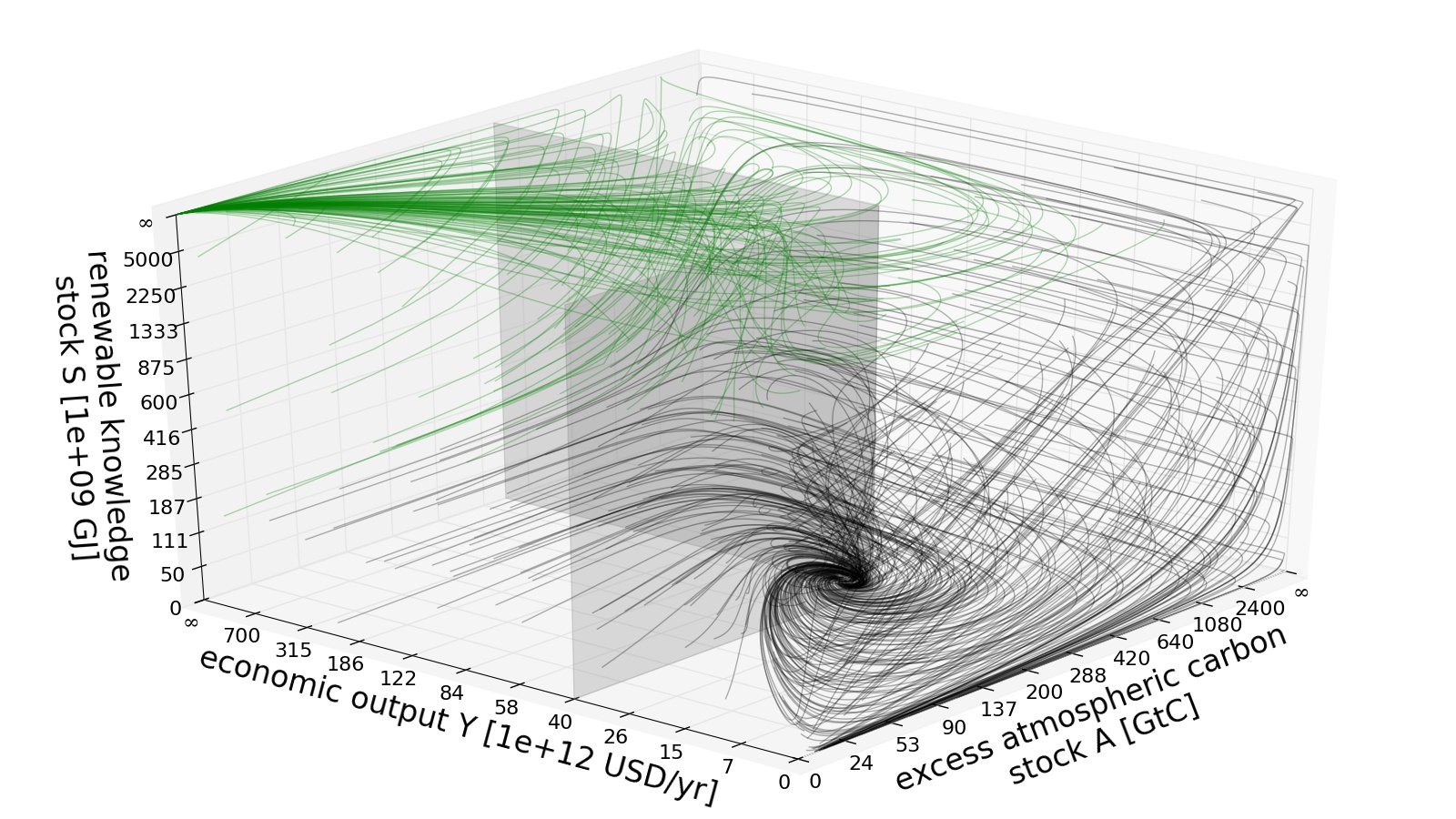}\label{fig:aws-default-flow}} \\
	\subfloat[Low growth management option ]{\includegraphics[width=0.8\columnwidth]{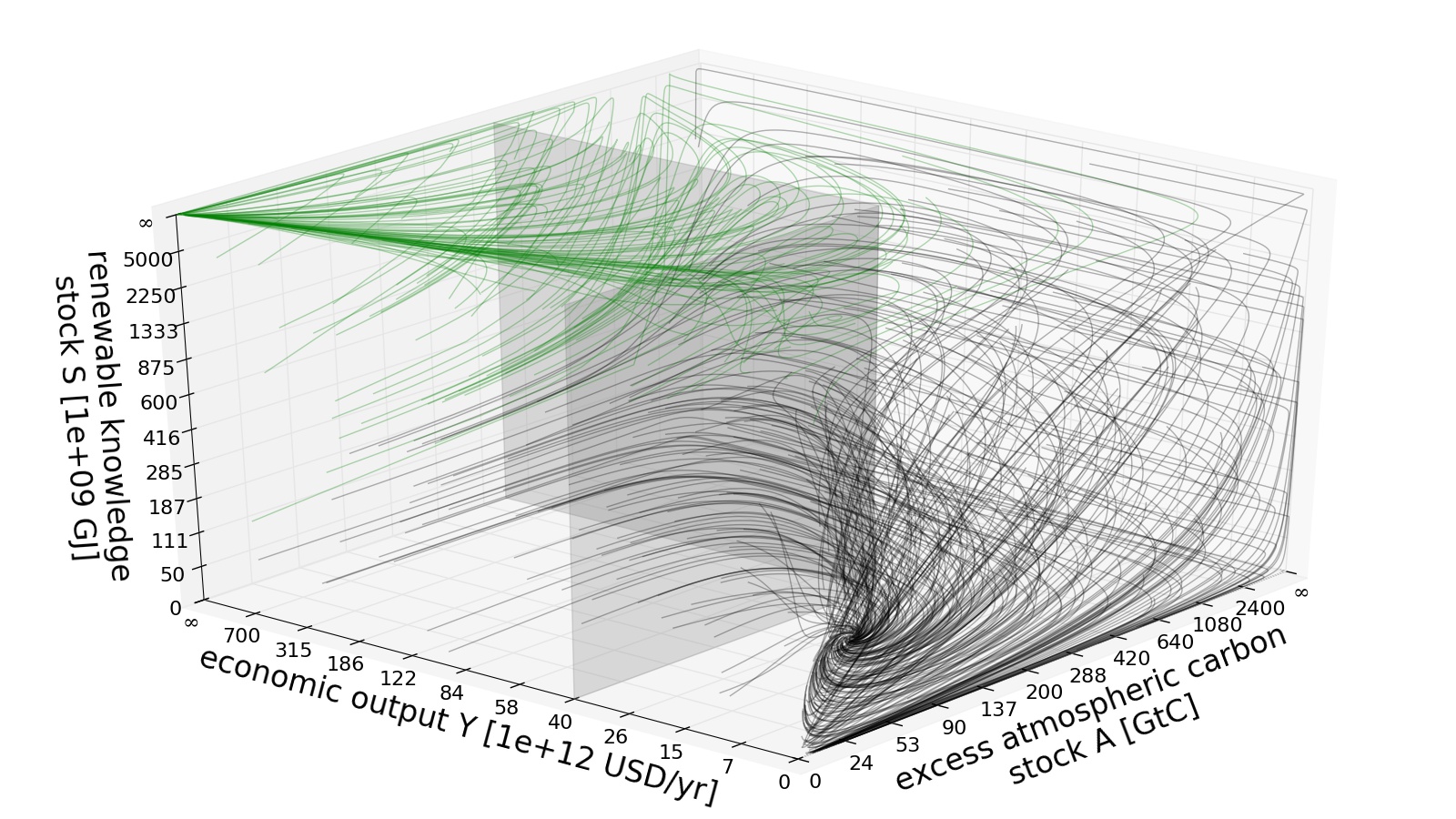}\label{fig:aws-degrowth-flow}} \\
	\subfloat[Energy transition management option ]{\includegraphics[width=0.8\columnwidth]{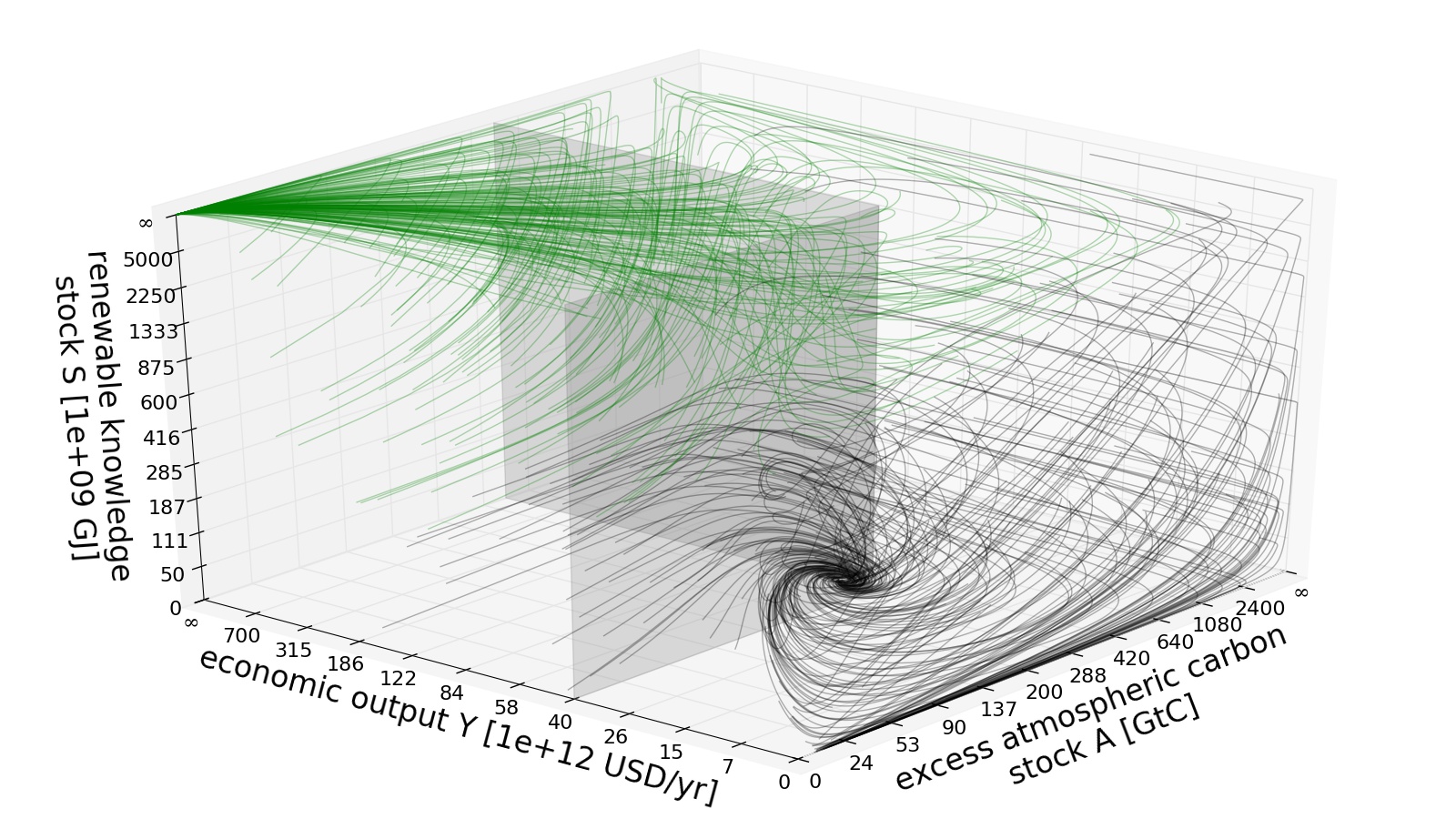}\label{fig:aws-et-flow}} \\
	\caption{The flows of the AYS model for (a) the default dynamics, (b) the low growth (\lg{}) and (c) the energy transition (\et{}) management options. The shaded planes indicate the two boundaries limiting the desirable region on the left of each panel. }
	\label{fig:aws-default+management}
\end{figure}


\subsection{Desirable states}
\label{sec:desirable-states}

Steffen~et~al.~\citep{Steffen2015} set the planetary boundary for climate change (\pbcc{}) to $350\,$ppm [ppm = parts per million] with an uncertainty zone until $450\,$ppm. We chose the desirable region to be where $A > A_{PB} = 345\,$GtC (above the pre-industrial level $A_0$), corresponding to the looser 450ppm  boundary (see \Cref{app:paramter-estimation} for the conversion).

The safe and just operating space by \citet{Raworth2012} complements the planetary boundaries with social foundations. Here we use a lower bound on the world gross product for illustration of a social foundation, even though we are well aware of the serious shortcomings of this indicator to measure decent livelihoods and basic human need satisfaction (e.g. \citep{Costanza2014,Jackson2016,Schmelzer2016}). Most importantly, the aggregate product does not capture any distributional aspects between countries and individuals. The exact value for such a lower boundary is open for discussion. For pragmatic reasons, we use here $Y_{SF }= 4\cdot 10^{13}\,$US\$ (\ysf{}), the economic production of the year 2000.

\Cref{fig:aws-default+management} shows the location of the boundaries in the three dimensional state space as two shaded planes. The desirable region is defined by the two boundaries and located on the left in all panels of \Cref{fig:aws-default+management}.


\subsection{Current state}
\label{sec:current-state}

The current state $\x_c = \left(A_c, Y_c, S_c\right)$ can be estimated as follows. $A_c$ is currently around 240\,GtC, corresponding to a concentration of 400\,ppm \citep{Betts2016}, and the world gross product of 2015 is around 70 Trillion US\$ \citep{WorldbankGDP}.
$S_c$ is estimated on the basis of the total past renewable energy consumption of roughly $2 \cdot 10^{12}$\,GJ \citep{Mongabay2015}.
This number has accumulated over roughly the same time as the characteristic knowledge depreciation time, $\tau_S = 50$\,a.
We assume roughly half of it has already depreciated, leaving $10^{12}$\,GJ.
Further accounting for the large error margins involved in estimating this number and because it contains hydroelectricity, whose growth potential is strongly limited, we take a conservative estimate by setting $S_c$ to a quarter of this values, $5 \cdot 10^{11}$\,GJ.

\begin{align}
\x_c = \left(\begin{array}{r}
A_c \\
Y_c \\
S_c
\end{array}\right)
= \left(\begin{array}{c}
240 \,\text{GtC} \\
7\cdot 10^{13} \,\text{US\$} \\
5\cdot 10^{11} \, \text{GJ}
\end{array}\right)
\end{align}

This estimation is quite rough and only provides limited information about the state of the target system because it reduces it to the three dimensions of a simplified model. However, it shows that the current state is still well inside the desirable region.


\subsection{Attractors}
\label{sec:aws-attractors}

With the above parameter values, the default dynamics of the AYS model has two fixed points.
The first fixed point has the coordinates at
\begin{align}
\x_b =
\begin{pmatrix} A_b \\ Y_b \\ S_b \end{pmatrix} =
\begin{pmatrix} \frac{\beta}{\theta} \\ \frac{\phi \epsilon \beta}{\theta \tau_A} \\ 0 \end{pmatrix} = 
\begin{pmatrix} 350\, \text{GtC} \\ 4.84\cdot 10^{13}\, \text{US\$} \\ 0 \end{pmatrix}.
\end{align}
It represents a fossil-based economy without renewable energy use, reduced economic output and constant climate damages. Therefore, we call it the ``black fixed point'' $\x_b$.

The second fixed point is located at the boundaries of the state space,
\begin{align}
\x_g = \begin{pmatrix} A_g \\ Y_g \\ S_g \end{pmatrix} = \begin{pmatrix} 0 \\ +\infty \\ +\infty \end{pmatrix},
\end{align}
with ``$+\infty$'' describing the boundary of the state space, as becomes clear in \Cref{sec:aws-compatification}.
This attractor corresponds to an economy with unbounded exponential growth of economic output and renewable knowledge. Also, fossil fuel usage and emissions decline towards zero. Therefore, we call this attractor the ``green fixed point''.

Both the black and green fixed points represent extreme asymptotic states. Even if they might correspond to unrealistic special cases of the modeled system, we find this acceptable because our focus lies on the model transients and the model only serves as an example.
With regard to the desirable regions, our analysis of the default dynamics shows that the green fixed point $\x_g$ is not violating either of the two boundaries, while the black one $\x_b$ violates the \pbcc{} as can be seen in \Cref{fig:aws-default-flow}.

Applying the \emph{low growth} management option moves the black fixed point to $\x_{b,\lg} = (A,Y,S)=(175\,$GtC$,2.42\cdot 10^{13}\,$US\$/a$,0)$, no longer violating the \pbcc{} (see \Cref{fig:aws-degrowth-flow}) but now violating the \ysf{}. Furthermore, for our choice of parameters and boundaries, no choice of the basic growth rate $\beta$ can make the black fixed point lie within the desirable region, thus satisfying both the \pbcc{}- and \ysf{}-boundaries at the same time.

Applying the \emph{energy transition} management option does not affect the location of the two attractors. But, more importantly, it changes the shape of the basins of attractions. When carefully inspecting \Cref{fig:aws-et-flow}, one can see that the volume the basin of attraction of the green fixed point is enlarged in comparison to the default flow in \Cref{fig:aws-default-flow}. Within the concept of Basin Stability \citep{Menck2013,Menck2014} the volume of the basin of attraction has been found to be an important indicator for an attractor's stability hence we will use a similar approach for the bifurcation analysis in \Cref{sec:analysis-bifurcation}.


\subsection{Dealing with the unbounded state space}
\label{sec:aws-compatification}

The model described in \Cref{sec:aws-model-description} and given by Eqs.~\eqref{eq:derivate-A}-\eqref{eq:definition-F-R-E} is defined on an unbounded state space. In this section, we map it to a bounded state space in order to deal with the diverging dynamics and apply the Saint-Pierre algorithm. We map the original state space $\X=\Rplus^3$ with the dynamical variables $\x = (A, Y, S)$ to a bounded space $\Y=[0,1)^3$ with transformed variables $\y = (a, y, s)$ and then add the point $\y_g = (0, 1, 1)$ which is the equivalent of $\x_g$ in the new coordinates. So we perform a change of coordinates
\begin{align}
\dot{\x} = f(\x) \quad \longrightarrow \quad \dot{\y} = F(\y),
\end{align}
where we switch from the old right-hand side (RHS) $f$ to the new RHS $F$. Following the explanations in \Cref{sec:unbounded-state-space-general}, we use the transformation
\bgroup
\renewcommand*{\arraystretch}{1.5}
\begin{align}
\begin{array}{@{}r@{\ }c@{}c@{\enspace}l@{}}
\Phi: & \X = \left[0, \infty\right)^3 &\longrightarrow& \Y = \left[0, 1\right)^3 \\
&A &\longmapsto& a = \frac{A}{A_{mid} + A}  \\
&Y &\longmapsto& y = \frac{Y}{Y_{mid} + Y}  \\
&S &\longmapsto& s = \frac{S}{S_{mid} + S},
\end{array}
\end{align}
\egroup
where the parameters $\x_{mid}=(A_{mid}, W_{mid}, S_{mid})^T$ are such that $\Phi(\x_{mid}) = (\frac{1}{2}, \frac{1}{2}, \frac{1}{2})^T$. They can be understood as the scale where the ``resolution is the best''. Hence, changing this value does not qualitatively influence the result, but a good choice can make them clearer.

This is exactly the transformation that has been used to create \Cref{fig:aws-default-flow-without-boundaries,fig:aws-default-flow,fig:aws-degrowth-flow,fig:aws-et-flow} and will be used for all the following figures, too. As we care most about the current state of the world, we choose $\x_{mid} = \x_c$.

Using \Cref{eq:general-coordinate-transformed-system}, we get a new set of ODEs with $\y = (a, y, s)$ as coordinates
\begin{subequations}
\begin{align}
\dot a &= \frac{ W_{mid}}{\phi \epsilon A_{mid}} \gamma (1-a)^2 \frac{y}{1-y} - \frac{a(1-a)}{ \tau_A}, \label{eq:rescaled-adot}\\
\dot y &= y(1-y) (\beta - \theta A_{mid} \frac{a}{1-a}),  \label{eq:rescaled-wdot} \\
\dot s &= (1-\gamma) \frac{W_{mid}}{\epsilon S_{mid}} (1-s)^2 \frac{y}{1-y} - \frac{s(1-s)}{\tau_S},  \label{eq:rescaled-sdot} \\[4mm]
\gamma &= \frac{(1-s)^\rho}{(1-s)^\rho + \left(\frac{S_{mid} s}{\sigma}\right)^\rho}, \label{eq:rescaled-gamma}
\end{align}
\end{subequations}
where $\gamma$ is the equivalent of $\Gamma$ in \Cref{eq:definition-gamma} but in the $\y$-coordinates.
The fixed points in the new $y$-coordinates are
\begin{subequations}
\begin{alignat}{3}
& \y_g = \begin{pmatrix} a_g \\ y_g \\ s_g \end{pmatrix} = \begin{pmatrix} 0 \\ 1 \\ 1 \end{pmatrix} \qquad && \longleftrightarrow \qquad && \x_g = \begin{pmatrix} A_g \\ Y_g \\ S_g \end{pmatrix} = \begin{pmatrix} 0 \\ \infty \\ \infty \end{pmatrix} \\
& \y_b = \begin{pmatrix} a_b \\ y_b \\ s_b \end{pmatrix} = \begin{pmatrix} \frac{\beta}{\beta + \theta A_{mid}} \\ \frac{\phi \epsilon \beta}{\phi \epsilon \beta + Y_{mid} \theta \tau_A} \\ 0 \end{pmatrix} \qquad && \longleftrightarrow \qquad && \x_b = \begin{pmatrix} A_b \\ Y_b \\ S_b \end{pmatrix} = \begin{pmatrix} \frac{\beta}{\theta} \\ \frac{\phi \epsilon \beta}{\theta \tau_A} \\ 0 \end{pmatrix}.
\end{alignat}
\end{subequations}
Now, we formally extend the dynamics such that $F(\y_g) = 0$.


\subsection{Results}
\label{sec:results-ays}

\begin{figure*}
\centering
\settrimmingRegions
\subfloat[]{\includegraphics[width=\columnwidth]{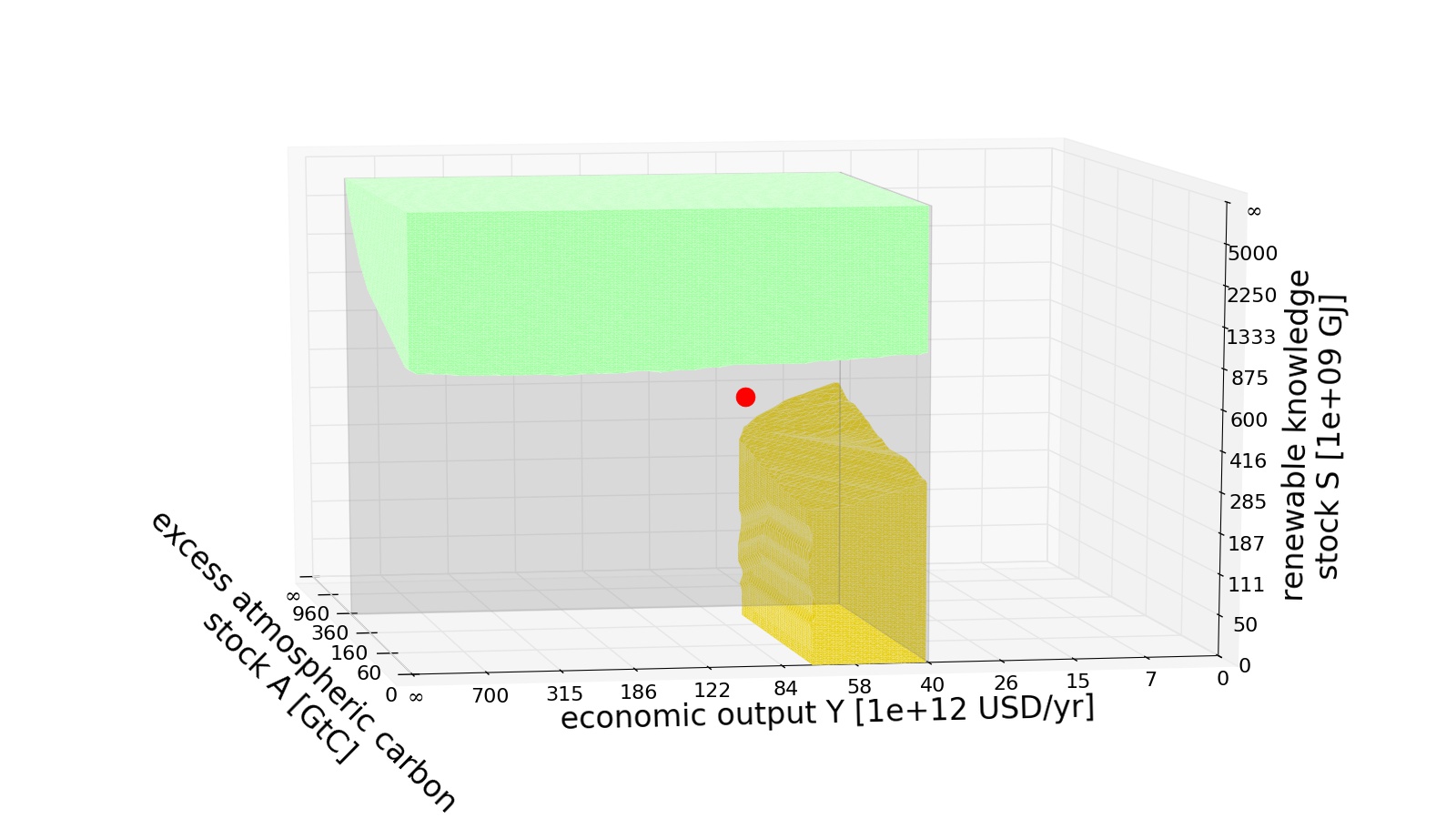}\label{fig:sc1-shelter-backwater}} \hfill
\subfloat[]{\includegraphics[width=\columnwidth]{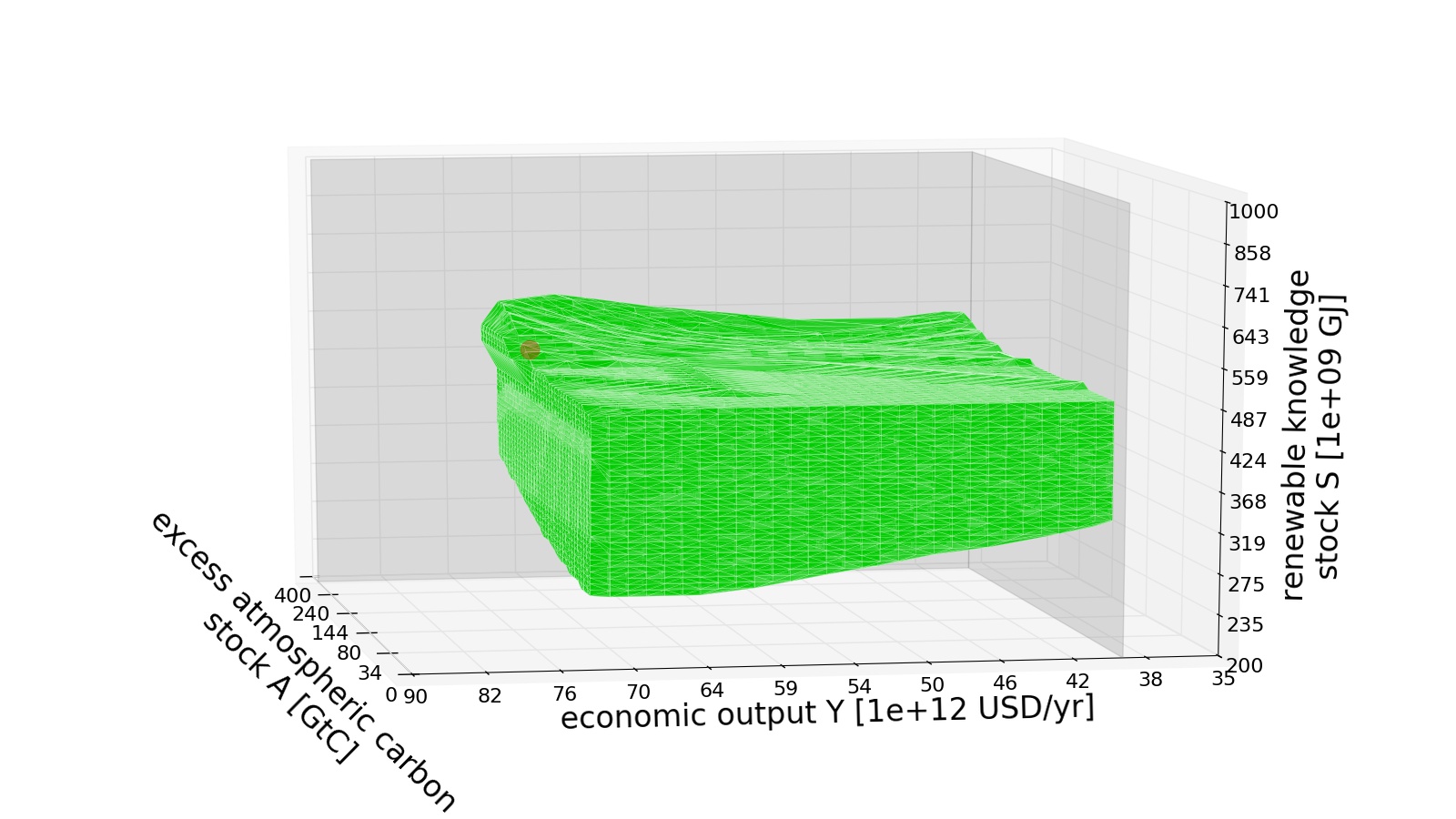}\label{fig:sc1-current-state-in-lake}}
\caption{Parts of the main \tsm{} partition for the AYS-model. The red dot in panel (a) indicates the estimation for the current state (see \Cref{sec:current-state}). In the shelter, the light green volume in panel (a), enough knowledge on renewable resource use has already been accumulated such that the relative price of renewable energy is low enough to satisfy the growing total energy demand without transgressing the boundaries. The \emph{low growth} management option allows the system to stay in the yellow backwater in panel (a) resulting in a \emph{zero-growth economy} within the desirable region. Panel (b) shows the time-limited lake (green) in a zoom into the state space around the current state.}
\label{fig:tsm-regions-ays-model}
\end{figure*}

We analyze the model in its compactified form \Cref{eq:rescaled-adot,eq:rescaled-wdot,eq:rescaled-sdot,eq:rescaled-gamma} using the nonlinear, local time-homogenization (\Cref{sec:time-homogenization}) and the Saint-Pierre algorithm (\Cref{sec:saint-pierre-algorithm}). Here, we do not provide the equations for the time-homogenized version because they are lengthy, their calculation straight-forward and they do not give much insight.
The most important identified regions are depicted in \Cref{fig:tsm-regions-ays-model}. We discuss their properties in the following.

The first regions to note are the shelter and the backwater shown in \Cref{fig:sc1-shelter-backwater}. 
In our model, the shelter, in which the system can stay in the desirable region without management, is the invariant kernel of the green fixed point's basin of attraction when restricting it to desirable states $\Xp$ only. 
States in this region correspond to a system that has accumulated enough knowledge for the energy production with renewable resources such that they become so cheap that there is no need for fossil fuels anymore. So the remaining (excess) CO$_2$ (above long-term equilibrium) is removed over time due to natural carbon uptake, leading the system to the green fixed point.
The model also contains a \emph{glade}, i.e., a region from which one can reach the shelter through the desirable region. It is a thin layer under the shelter and has not been visualized in \Cref{fig:tsm-regions-ays-model}.

The \emph{backwater}, where one can stay in the desirable region forever but needs to apply management over and over again, is the partition of the desirable region where the growth of economic output and hence of emissions can be restricted.
That way, the atmospheric carbon concentration can be kept within the planetary boundary.
Within the backwater, decarbonization of the economy is impossible because the given maximal carbon tax or renewable subsidy policies are too weak to make renewables competitive with fossil fuel.
Instead, one can use the low growth management option to stay in a state that corresponds to a carbon-based economy, in which climate impacts are compensated by low economic growth. The atmospheric carbon level is relatively high but still within the boundary and in equilibrium with the emissions.
Using the low growth option, one can therefore stay within the desirable region but cannot reach the green fixed point.

For simplicity, we did not include the option to choose a value of the base growth rate lying between the two options $\beta$ and $\beta_\lg$. So formally, the management strategy required to stay in the desirable region involves switching between $\beta$ and $\beta_\lg$ because either of these two values alone leads to a black fixed point in the undesirable region.
Still, it is easy to see that this management strategy is equivalent to using a constant intermediate value of $\beta_m$ (e.g.\ $\beta_m=2.7\,\%/$a) instead.
The dynamics with such an intermediate growth rate has a black fixed point that lies in the desirable region. So this option is implicitly included in the model.
In other words, one only needs to model the maximal and minimal values of the option space explicitly because the \tsm{}-framework allows for arbitrarily fast switches between the management options. Intermediate options, more precisely all convex combinations of options, are implicitly taken into account.
Still, replacing the value of $\beta_\lg$ with $\beta_m$ would introduce further changes. The maximal management in the transient would be restricted and thus the size of the backwaters reduced. We show this in \Cref{sec:analysis-bifurcation}. Here, it becomes obvious that the \tsm{}-partition depends both on the asymptotics and the transient dynamics of the model.

According to our analysis, the current state $\x_c$ (as estimated in \Cref{sec:current-state}) lies between the shelter and the backwater discussed above, in a \emph{time-limited lake} as the zoom into the state space in \Cref{fig:sc1-current-state-in-lake} shows. Therefore, the results suggest that humanity currently faces a so-called lake dilemma in the AYS model. In such a dilemma, we have to make the qualitative decision between staying within the desirable region uninterruptedly but depending on management forever or going through the undesirable region to finally end up in the shelter. In the model, the choice is between using the energy transition option in order to speed up the knowledge accumulation on renewable resource use and finally reach the green fixed point or the low growth option that restricts the total energy demand and helps reaching the black fixed point. Even with a combined usage of the different options, the model does not allow to reach the shelter without transgressing the boundaries when going for the first option.
This classification of $\x_c$ depends on the finite resolution $h$ of the applied algorithm. An alternative analysis of the AYS model with deep reinforcement learning suggest that $\x_c$ is in a glade region instead of a lake region \cite{strnad2020}.


\subsection{Bifurcation analysis}
\label{sec:analysis-bifurcation}

\begin{figure}
\centering 
\subfloat[Growth control]{\includegraphics[width=0.4\columnwidth, valign=c]{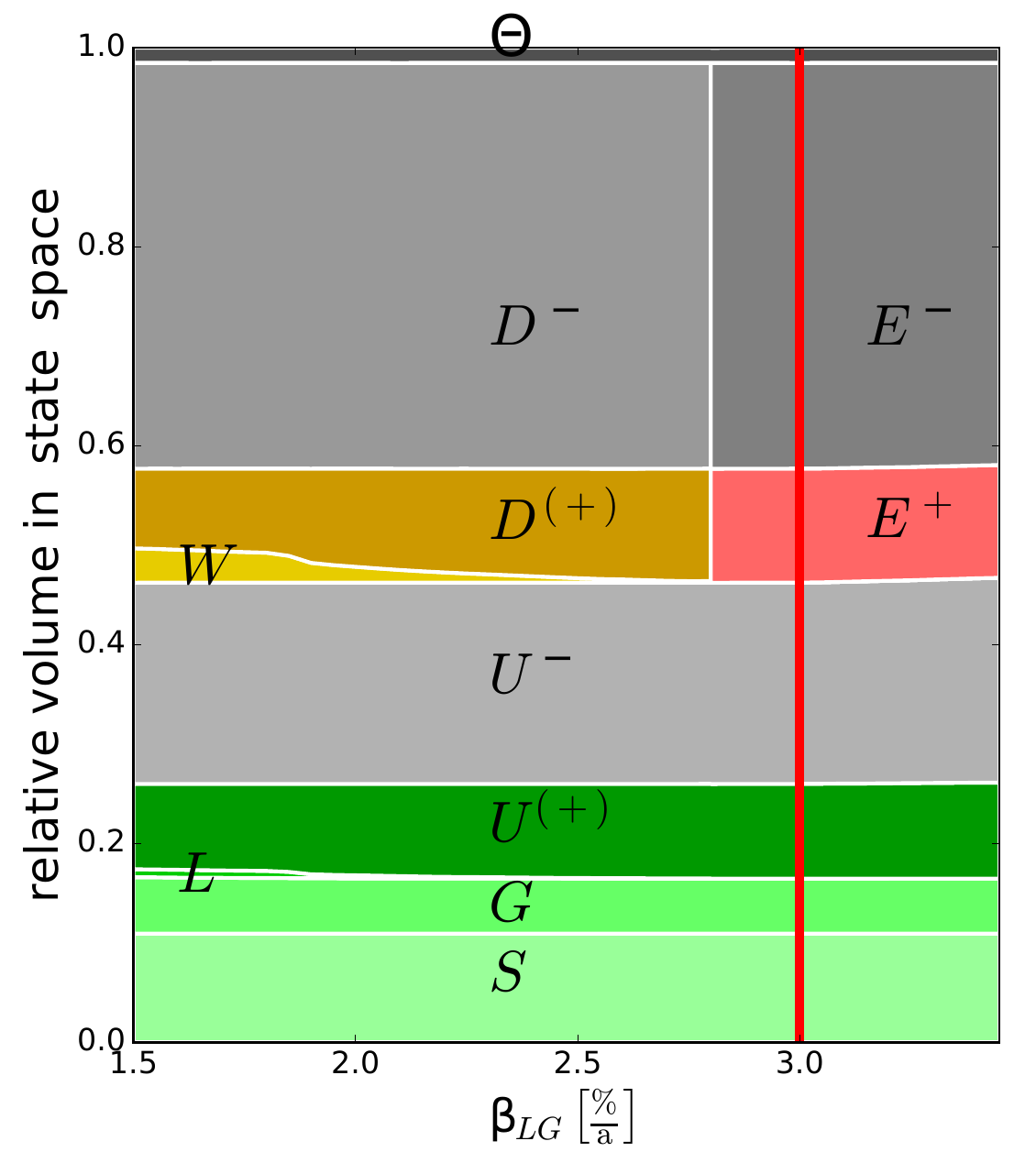}\label{fig:sc1-bifurcation-gr}}
\hfill
\subfloat[Energy transition control]{\includegraphics[width=0.4\columnwidth, valign=c]{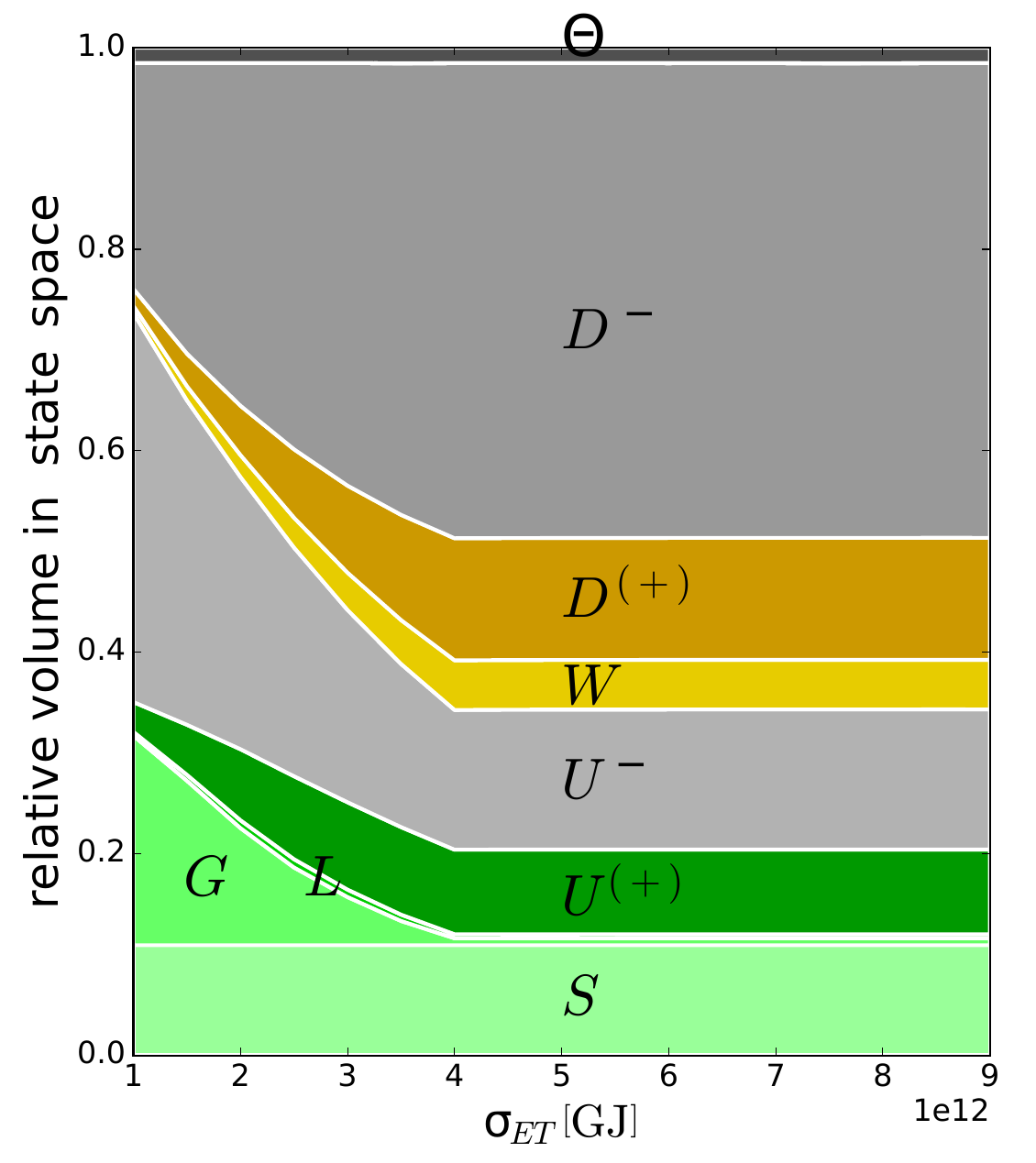}\label{fig:sc1-bifurcation-et}}
\hfill
\subfloat{\includegraphics[width=0.19\columnwidth, valign=c]{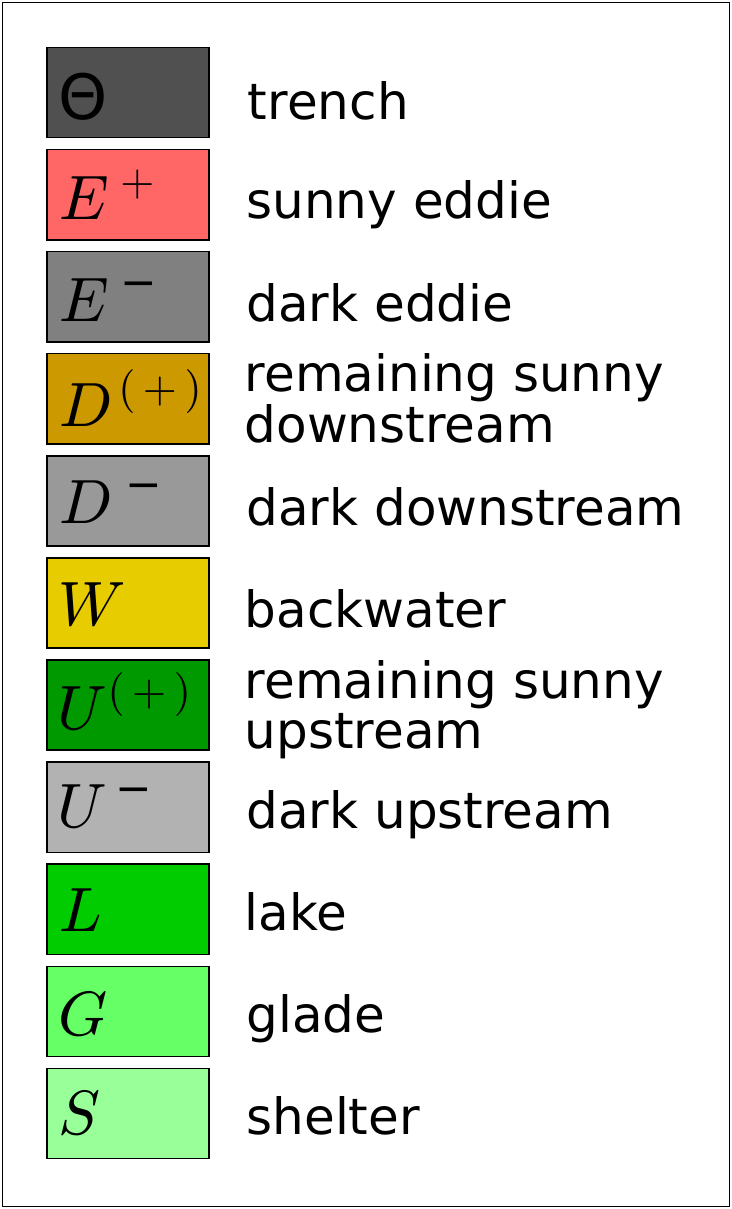}}
\caption{The bifurcation diagrams show the relative volume of different regions of the \tsm{} partition over the parameters of the two management options.
Panel (a) shows a \emph{downstream-eddies bifurcation} when varying the manageable growth rate $\beta_\lg$. The default growth rate is marked by a red vertical line.
Panel (b) indicates strong increases in the size of the glade and dark upstream with stronger management of the energy transition, corresponding to lower relative costs for renewable energy and smaller values of $\sigma_{ET}$.
}
\end{figure}

Varying the parameters corresponding to the two management options may lead to what can be understood as {\em \tsm{} bifurcations} with respect to the topology of \tsm{} regions in state space.
A full theory of \tsm{} bifurcations would be based on a suitable notion of topological equivalence between control systems with default trajectories and specified desirable regions, just like the standard theory of bifurcations is based on such an equivalence between (uncontrolled) dynamical systems \cite{Kuznetsov1998}. Because this is beyond the scope of this paper, we use a working definition of \tsm{} bifurcation here that is based on sudden, qualitative changes of the \tsm{} partition at certain critical parameters (see \Cref{app:tsm-bifurcation-definition} for details).

This allows us to use the relative volume of each region an indicator for \tsm{} bifurcations, motivated by the concept of Basin Stability \citep{Menck2013,Menck2014} and its extensions \citep{Hellmann2015,Klinshov2015,Kittel2017,Mitra2015,Van2016,Mitra2017,Mitra2017a}.
We use uniformly distributed points in state space for the Saint-Pierre algorithm. Therefore, we can estimate the relative volume of one region with the number of points associated to this region over the total number.

When varying $\beta_\lg$ corresponding to the \emph{low growth} management option (\lg{}) from $1.5\,\%/$a to $3.5\,\%/$a, a \emph{downstream-eddies bifurcation} occurs. Until the fixed point of the \lg{} flow crosses the planetary boundary at the critical value of $\beta_\lg =$ $ \beta_{PB}\approx $ $ 2.95\%/$a there exists a \emph{backwater}\footnote{Because of the discretization of state space and time in the estimation, \Cref{fig:sc1-bifurcation-gr} suggests that the bifurcation occurs already at $\approx 2.8\%$/a. However, this is a numerical inaccuracy.}. Beyond this point, there are only \emph{eddies} left because the focus of the default flow and the one of the low-growth flow are both in the undesirable region. Thus, for $\beta_\lg > \beta_{PB}$ it is only possible to switch between the two flows such that the system circles around both foci and visits the desirable region again and again, while having to pass through the undesirable region in between.

As discussed in \Cref{sec:results-ays}, the \emph{backwater} can be explained by the fact that there is an in-between value $\beta_m = 2.7\,\%/$a, for which the black fixed point lies in the desirable region. \Cref{fig:sc1-bifurcation-gr} shows that the volume of the backwater for this intermediate value $\beta_m = 2.7\,\%/$a is smaller than for $\beta_\lg = 1.5\,\%/$a. This is due to the restriction of the maximally available management such that it becomes impossible to stay within the desirable region forever for a large share of trajectories starting of in the downstream.

\Cref{fig:sc1-bifurcation-et} depicts the changes in relative volume of the \tsm{}-partition over the energy transition control $\sigma_{ET}$. Especially for a strong management corresponding to small values of $\sigma_{ET}$, we observe strong changes in the relative volume of the glade, the lake and the backwater. The glade increases in size, because it is possible to reach the green fixed point from more initial conditions without transgressing the boundaries. 

\section{Summary \& Outlook}
\label{sec:summary}

In this paper, we introduce a variant definition of the topology of sustainable management (\tsm{}) in terms of viability theory that operationalizes the partition of state space into qualitatively different regions. We apply the Saint-Pierre algorithm to estimate the volume and visualize the position and shape of different \tsm{} regions.
Because the algorithm works on bounded spaces only, we use a coordinate transformation to a bounded space that best resolves a fixed scale of interest. Furthermore, we address the problem that time scales vary by orders of magnitude by introducing a nonlinear, local time homogenization.

We apply these novel concepts to a three dimensional example model of interactions between climate change, economic growth and energy transition.
Even though we design the model in a stylized way, the analysis reveals a rich \tsm{} topology.
This results from the two fixed points in the model corresponding to a fossil and a green economy.
We estimate the current state of the world to lie inside a finite-time lake. In the model, humanity is therefore facing a pressing lake-dilemma: it has limited time to make a choice between two qualitatively different options.
Moreover, we analyze how the relative volume of different \tsm{} regions change with the management parameters. We identify a downstream-eddies bifurcation, in which a downstream region becomes an eddie after one of the fixed points crosses the boundary of the desirable region.

The example model presented here has a comparatively low-complexity that allows us to focus on the application of the novel methods for estimating the \tsm{} regions. Relevant processes, like the carbon cycle, economic dynamics and different energy productions, are represented only in a stylized way. Future work can apply our framework to extended versions of our model. For example, it could incorporate other management options, such as carbon capture and storage, or introduce higher dimensional representations of the carbon cycle or economic system. The models developed in Refs. \citep{Nitzbon2017, Nitzbon2016thesis} could to be good candidates for such extensions because they have only a limited number of dynamical variables.

Extending the framework also to high-dimensional models would require improved algorithms for the estimation of viability kernels and capture basins.
While there is work on efficient algorithms for such problems \citep{Bokanowski2006,Maidens2013,Brias2016,Alvarez2016}, they need to be adjusted and extended to be applicable for the computation of the \tsm{} partition.
Furthermore, future theoretical research could further formalize our working definition for \tsm{} bifurcations and extend the \tsm{} framework to stochastic dynamical systems \citep{Anishchenko2007}.
Such further developments of the presented methods would be highly valuable to advance our understanding of complex social-ecological models of the Earth system \citep{Mueller-Hansen2017,Donges2020} that can guide sustainable development within a safe and just operating space.

\subsection*{Software availability}

For this project, an open source Python-library for viability computations called \emph{pyviability} has been developed and is available online at \href{https://timkittel.github.io/PyViability/}{https://timkittel.github.io/PyViability/} (published July 7, 2017). The source code for the model-specific computations including a description to reproduce the results presented in this article is available at \href{https://github.com/timkittel/ays-model/}{https://github.com/timkittel/ays-model/} (published July 7, 2017).

\begin{acknowledgement}

\subsection*{Acknowledgement}

This article was developed within the scope of the IRTG 1740/TRP 2011/50151-0, funded by the DFG/FAPESP.
This work was conducted in the framework of PIK’s flagship project on coevolutionary pathways (\textsc{copan}).
The authors thank CoNDyNet (FKZ 03SF0472A) for their cooperation.
The authors gratefully acknowledge the European Regional Development Fund (ERDF), the German Federal Ministry of Education and Research and the Land Brandenburg for supporting this project by providing resources on the high performance computer system at the Potsdam Institute for Climate Impact Research.
The authors thank the developers of the used software: Python \citep{python}, Numerical Python \citep{numpy} and Scientific Python \citep{scipy}.

The authors thank
Sabine Auer,
Wolfram Barfuss,
Karsten Bölts,
Catrin Ciemer,
Eduardo Costa,
Jonathan Donges,
Reik Donner,
Jaap Eldering,
Jasper Franke,
Roberto Gueleri,
Frank Hellmann,
Jakob Kolb,
Julien Korinman,
Till Koster,
Chiranjit Mitra,
Jan Nitzbon,
Ilona Otto,
Tiago Pereira,
Camille Poignard,
Francisco Rodrigues,
Edmilson Roque dos Santos,
Stefan Ruschel,
Paul Schultz,
Fabiano Berardo de Sousa,
Lyubov Tupikina,
and Kilian Zimmerer
for helpful discussions and comments.
\end{acknowledgement}

\printbibliography



\appendix
\renewcommand*{\thesection}{\Alph{section}} 
{
\section*{Appendix}
\addcontentsline{toc}{section}{Appendix}  
}

\section{Existence of eddies}
\label{app:existence-eddies}

\paragraph{Definition of eddie-like sets:} We call a pair of sets $\A^{+/-} \subseteq \X$ \emph{eddie-like} if and only if they fulfill the following two conditions: (i) $\A^{+/-} \subseteq \Xpm - \topU - \topD$ and (ii) $\A^{+/-} \subseteq \text{Capt}\left(A^{-/+}\right)$. Note the inverted order of the signs in the last term.

\paragraph{Union of two eddie-like pair sets are also eddie-like:}
For two eddie-like pairs of sets $\topEpm_1$ and $\topEpm_2$, the union pair $\topEpm_3 = \topEpm_1 \cup \topEpm_2$ is eddie-like, too.

\paragraph{Proof:} The first condition is trivially fulfilled and the second one follows straight away from $\text{Capt}\left(\A\right) \cup \text{Capt}\left(\B\right) = \text{Capt}\left(\A \cup \B\right)$ for two state sets $\A, \B \subseteq \X$.
Hence, the union of all eddie-like pairs of sets is maximal and \emph{eddies} exist. \qed

\section{Working definition of \tsm{} bifurcations}
\label{app:tsm-bifurcation-definition}

For the analysis in \Cref{sec:analysis-bifurcation}, we use the following working definition of \tsm{} bifurcations.
Let $\psi=(\psi_1,\dots,\psi_k)$ be a vector of parameters of the control system, let $\topS(\psi),\topG(\psi),\dots,\topT(\psi)\subseteq\X$ be the various \tsm{} regions in dependence on parameter values, and let $f(\psi)=(\topS(\psi),\topG(\psi),\dots,\topT(\psi))$ be the function that maps parameter values to the tuple of these regions.
Now, we say that at a certain point $\psi^\ast$ in parameter space a {\em \tsm{} bifurcation} happens if and only if $f(\psi)$ is discontinuous at $\psi^\ast$ with respect to the Hausdorff metric between subsets of $\X$.
In particular, we have a \tsm{} bifurcation if the volume of some of the regions changes discontinuously at $\psi^\ast$. 

\section{Parameter estimation}
\label{app:paramter-estimation}

To get a roughly realistic setting, we estimated the parameters of the model using several publicly available data sources.

$A_0$ was taken from \citep{Ciais2013} and slightly rounded. 
$\tau_A$ and $\beta$ were taken from \citep{Kellie-Smith2011a}.
$\phi$ was based on the ton oil equivalent of various fossil fuels and a typical mass share of 90\% carbon in fossil fuels, as described in \citep{Nitzbon2016thesis}.

Assuming that two degrees warming correspond to a carbon concentration of around 450\,ppm \citep{IPCC2013}
and thus to a carbon stock of 950\,GtC (both being 1.6 times their pre-industrial value),
we require that the total growth rate $\beta - \theta A_1$ becomes zero 
for $A_1 \approx 950\,$GtC${} - A_0 = 350\,$GtC,
hence $\theta$ was taken to be $\beta / A_1$ $\approx $ $8.57\cdot 10^{-5}/{}$(GtC a).

$\epsilon$ was estimated from the World Bank's primary energy intensity data \citep{WorldbankEnergy}.

For $\tau_S$, the characteristic depreciation time of renewable energy knowledge, no reliable source was found, 
so we made a very coarse guess by setting it roughly to the length of an average working life of 50\,a.

The break-even knowledge level $\sigma$ was also estimated very coarsely.
According to \citep{Mongabay2015},
past cumulative world consumption of renewable energy is $\approx 2\cdot 10^{18}\,$Btu $\approx 2\cdot 10^{12}\,$GJ
or roughly 20 years of world energy consumption.
To be on the conservative side and avoid overestimating the potential of renewables,
we took $\sigma$ to be two times that value.

$\rho$ was set as follows. 
We assume fossil and renewable energy production costs of 
$C_F\propto F^{1+\gamma}$ and $C_R\propto R^{1+\gamma} / S^\lambda$,
where $\gamma > 0$ is a convexity parameter and $\lambda > 0$ is a learning exponent.
Then energy prices are
$\pi_F\propto\partial C_F/\partial F$ $\propto F^\gamma$
and
$\pi_R\propto\partial C_R/\partial R$ $\propto R^\gamma / S^\lambda$.
In the price equilibrium, $\pi_F = \pi_R$, 
hence $R / F \propto S^{\lambda / \gamma}$,
and thus $\rho = \lambda / \gamma$.
According to \citep{Rubin2015}, the learning rate $LR = 1 - 2^{-\lambda}$
of several renewables is around 1/8, 
hence $\lambda \approx \log_2 (8/7) \approx 0.2$.
Assuming a mild convexity of $\gamma\approx 0.1$, we get $\rho\approx 2$.


\end{document}